\documentclass[12pt,oneside,reqno]{amsart}
\usepackage{indentfirst}
\usepackage{graphicx}
\usepackage{mathrsfs}
\usepackage{stmaryrd}
\usepackage{amsfonts}
\usepackage{enumerate,amsmath,amssymb,amsthm}
\usepackage{color}
\usepackage{float}
\newcommand{\dif}{\mathrm{d}}

\newcommand{\be}{\begin{eqnarray}}
\newcommand{\ee}{\end{eqnarray}}
\newcommand{\ce}{\begin{eqnarray*}}
\newcommand{\de}{\end{eqnarray*}}
\newtheorem{theorem}{Theorem}[section]
\newtheorem{lemma}[theorem]{Lemma}
\newtheorem{remark}[theorem]{Remark}
\newtheorem{definition}[theorem]{Definition}
\newtheorem{proposition}[theorem]{Proposition}
\newtheorem{Examples}[theorem]{Examples}
\newtheorem{corollary}[theorem]{Corollary}
\def\e{\varepsilon}
\def\s{\sigma}
\def\t{\theta}
\def\a{\alpha}

\def\d{\delta}
\def\p{\partial}
\def\g{\gamma}

\def\[{{\Big[}}
\def\]{{\Big]}}
\def\<{{\langle}}
\def\>{{\rangle}}
\def\({{\Big(}}
\def\){{\Big)}}

\def\no{\nonumber}
\def\bt{\begin{theorem}}
\def\et{\end{theorem}}
\def\bl{\begin{lemma}}
\def\el{\end{lemma}}
\def\br{\begin{remark}}
\def\er{\end{remark}}
\def\bx{\begin{Examples}}
\def\ex{\end{Examples}}
\def\bd{\begin{definition}}
\def\ed{\end{definition}}
\def\bp{\begin{proposition}}
\def\ep{\end{proposition}}
\def\bc{\begin{corollary}}
\def\ec{\end{corollary}}

\def\cD{{\mathcal D}}

\def\mE{{\mathbb E}}
\def\mF{{\mathbb F}}

\def\mN{{\mathbb N}}

\def\mP{{\mathbb P}}

\def\mR{{\mathbb R}}

\def\sA{{\mathscr A}}

\def\sF{{\mathscr F}}

\def\sL{{\mathscr L}}

\def\sV{{\mathscr V}}

\def\geq{\geqslant}
\def\leq{\leqslant}
\allowdisplaybreaks[3]

\begin{document}
	
\allowdisplaybreaks
	
\title{Average principles for forward-backward multivalued stochastic systems and homogenization for systems of nonlinear parabolic PDEs }
	
\author{Huijie Qiao}
	
	
	
\thanks{This work was supported by NSF of China (No.12071071)}
		
	
\date{}
\dedicatory{School of Mathematics,
		Southeast University\\
		Nanjing, Jiangsu 211189, China\\
		hjqiaogean@seu.edu.cn}

\begin{abstract}
This work concerns about forward-backward multivalued stochastic systems. First of all, we prove one average principle for general stochastic differential equations in the $L^{2p}$ ($p\geq 1$) sense. Moreover, for $p=1$ a convergence rate is presented. Then combining general stochastic differential equations with backward stochastic variation inequalities, we establish the other average principle for backward stochastic variation inequalities in the $L^{2}$ sense through a time discretization method. Finally, we apply our result to nonlinear parabolic partial differential equations and obtain the homogenization of them.  
\end{abstract}

\maketitle \rm

\section{Introduction}

Given a complete probability space $(\Omega, \sF, \mP)$ on which a $l$-dimensional standard Brownian motion $W$ is defined. Assume that $\mF:=(\sF_t)_{t\geq 0}$ is the $\mP$-augmentation of the natural filtration of $W$. Fix $T>0$ and consider the following SDE: $0\leq t\leq s\leq T$
\be\left\{\begin{array}{l}
\dif X_{s}^{\e,t,\zeta}=b(\frac{s}{\e},X_{s}^{\e,t,\zeta})\dif s+\s(\frac{s}{\e},X_{s}^{\e,t,\zeta})\dif W_{s},\\
X_{t}^{\e,t,\zeta}=\zeta,
\end{array}
\right.
\label{multfsde}
\ee
where these mappings $b: \mR_+\times\mR^{m}\rightarrow\mR^{m}$ and $\s: \mR_+\times\mR^{m}\rightarrow\mR^{m\times l}$ are all Borel measurable. $\zeta$ is a $\sF_t$-measurable random variable and $0<\e<1$ is a small parameter. Eq.(\ref{multfsde}) is usually called a multiscale stochastic differential equation (SDE for short). People care about the limit of its solution as $\e\rightarrow 0$, that is, the average principle. Nowadays, there have been many results about average principles about Eq.(\ref{multfsde}). For example, in \cite{nn} N'Goran and N'Zi studied the average principle for multivalued SDEs in the probability sense. Later, Xu and Liu \cite{xl} improved the result in \cite{nn} to the convergence in the $L^2$ sense. Recently, if $b, \s$ satisfy the local Lipschitz and monotone conditions, Guo et al. \cite{gxwh} established an average principle for Eq.(\ref{multfsde}) in the $L^2$ sense. And Shen, Song and Wu \cite{ssw} also proved an average principle in the $L^2$ sense, when $b, \s$ depend on the distribution of $X_{s}^{\e,t,\zeta}$. Here we show an average principle for Eq.(\ref{multfsde}) in the $L^{2p}$ sense. Moreover, for $p=1$ a convergence rate is presented. 

Next, we couple Eq.(\ref{multfsde}) with a backward stochastic variation inequality and investigate the limit for the forward-backward multivalued stochastic system as $\e\rightarrow 0$. 
Concretely speaking, consider the following backward stochastic variation inequality (SVI for short): $0\leq t\leq s\leq T$
\be\left\{\begin{array}{l}
\dif Y_{s}^{\e,t,\zeta}\in \p\varphi(Y_{s}^{\e,t,\zeta})\dif s-[f_1(\frac{s}{\e},X_{s}^{\e,t,\zeta},Y_{s}^{\e,t,\zeta})+f_2(Z_{s}^{\e,t,\zeta})]\dif s+Z_{s}^{\e,t,\zeta}\dif W_{s},\\
Y_{T}^{\e,t,\zeta}=g(X_{T}^{\e,t,\zeta})\in\overline{\cD(\p\varphi)},
\end{array}
\right.
\label{multfbsde}
\ee
where $\varphi$ is a lower semicontinuous convex function, $\p \varphi$ is its subdifferential operator, and these mappings $f_1: \mR_+\times\mR^{m}\times\mR^{d}\rightarrow\mR^{d}$, $f_2: \mR^{d\times l}\rightarrow\mR^{d}$ and $g: \mR^m\rightarrow\mR^d$ are all Borel measurable. If $f_1(s,x,y)$ is independent of $s$ and $f_2(z)=0$, Essaky and Ouknine \cite{eo} concluded that if the solution $X^{\e,t,\zeta}$ of Eq.(\ref{multfsde}) converges in law to the solution of the corresponding average equation, the solution $Y^{\e,t,\zeta}$ of Eq.(\ref{multfbsde})  also converges in law to the solution of the corresponding average equation. Recently, Hu, Jiang and Wang \cite{hjw} studied Eq.(\ref{multfbsde}) with $\varphi=0$. However, there they didn't explicitly express the convergence of $X^{\e,t,\zeta}$ and $Y^{\e,t,\zeta}$. Here we prove that $Y^{\e,t,\zeta}$ converges to the solution of the corresponding average equation in the $L^2$ sense. Besides, we mention that if $f_1(s,x,y)$ is independent of $x$ and $g(X_{T}^{\e,t,\zeta})$ is replaced by a random variable, a convergence rate can be obtained. But here we don't give out it. This is because we want to apply the average principle for Eq.(\ref{multfbsde}) to the homogenization of nonlinear parabolic partial differential equations (PDEs for short).  

In the following, we consider the following nonlinear PDE:
\be\left\{\begin{array}{l}
\frac{\p u^\e(t,x)}{\p t}+\sL^\e u^\e(t,x)+f_1(\frac{t}{\e},x,u^\e(t,x))+f_2(\nabla u^\e(t,x)\s(\frac{t}{\e},x))\in \p \varphi(u^\e(t,x)), ~t\in[0,T], \\ 
u^\e(T,x)=g(x), \quad u^\e(t,x)\in\overline{{\rm Dom}(\varphi)}, \quad x\in\mR^m,
\end{array}
\right.
\label{1dpdein1}
\ee
where 
\ce
\sL^\e:=\frac{1}{2}\sum\limits_{i,j=1}^m (\s\s^*)_{ij}(\frac{t}{\e},x)\frac{\p^2}{\p x_i \p x_j}+\sum\limits_{i=1}^m b_i(\frac{t}{\e},x)\frac{\p}{\p x_i}.
\de
Eq.(\ref{1dpdein1}) is often used to model some obstacle and control problems (cf. \cite{eo}). If $\varphi=0$ and $b,\s, f_1$ are independent of $s$, there have been many results about the homogenization of Eq.(\ref{1dpdein1}) (cf. \cite{blp, bh, bhp, bi, d, p}). Let us recall some works. In \cite{blp}, Bensoussan, Lions and Papanicolaou elaborated systematically the homogenization for PDEs with periodic structures by analytical methods. Later, Buckdahn and Ichihara \cite{bi} studied the homogenization for Hamilton-Jacobi-Bellman equations with periodic structures by a probability approach. If $\varphi=0$ and $b,\s, f_1$ depend on $s$, Hu, Jiang and Wang \cite{hjw} proved the homogenization for Eq.(\ref{1dpdein1}) by a probability approach. If $\varphi\neq 0$, $f_1(s,x,y)$ is independent of $s$ and $f_2(z)=0$, Essaky and Ouknine \cite{eo} observed the homogenization for Eq.(\ref{1dpdein1}) by a probability approach. However, up to now, there are few results about the homogenization of Eq.(\ref{1dpdein1}). Here we investigate the homogenization for Eq.(\ref{1dpdein1}) and the other two types of PDEs. 

The novelty of this paper lies in three folds. The first fold is that we prove one average principle for general stochastic differential equations in the $L^{2p}$ ($p\geq 1$) sense. Moreover, for $p=1$ a convergence rate is presented, which is important for numerical simulation. The second fold is that we establish the average principle for backward stochastic variation inequalities in the $L^{2}$ sense through a time discretization method. Since the operator $\p \varphi$ is multivalued, nonlinear and not smooth, we need some new ideas and approaches so as to reach this goal. The third fold is that we obtain the homogenization of nonlinear parabolic PDEs without periodic structures. Therefore, our result is more general.

Finally, this paper is arranged as follows. In the next section, we introduce notations and concepts, and recall some results used in the sequel. In Section \ref{ap}, two main results are formulated. Then the proofs of main results are placed in Section \ref{xbarxdiffproo} and \ref{averprinproo}. In Section \ref{app}, we apply our result to nonlinear parabolic PDEs. Finally, we give an example to explain our result in Section \ref{exam}.
 
The following convention will be used throughout the paper: $C$ with or without indices will denote different positive constants whose values may change from one place to another.

\section{Preliminary}\label{fram}

In this section, we introduce notations and concepts, and recall some results used in the sequel.

\subsection{Notations}\label{nota}

In this subsection, we introduce some notations.

For convenience, we shall use $\mid\cdot\mid$ and $\parallel\cdot\parallel$  for norms of vectors and matrices, respectively. Furthermore, let $\langle\cdot$ , $\cdot\rangle$ denote the scalar product in $\mR^d$. Let $B^*$ denote the transpose of a matrix $B$.

Let $C(\mR^d)$ be the collection of continuous functions on $\mR^d$ and $C^2(\mR^d)$ be the space of continuous functions on $\mR^d$ which have continuous partial derivatives of order up to $2$.

\subsection{Maximal monotone operators}\label{mmo}

In this subsection, we introduce maximal monotone operators. 

For a multivalued operator $A: \mR^d\mapsto 2^{\mR^d}$, where $2^{\mR^d}$ stands for all the subsets of $\mR^d$, set
\ce
&&\cD(A):= \left\{x\in \mR^d: A(x) \ne \emptyset\right\},\\
&&Gr(A):= \left\{(x,y)\in \mR^{2d}:x \in \cD(A), ~ y\in A(x)\right\}.
\de
We say that $A$ is monotone if $\langle x_1 - x_2, y_1 - y_2 \rangle \geq 0$ for any $(x_1,y_1), (x_2,y_2) \in Gr(A)$, and $A$ is maximal monotone if 
$$
(x_1,y_1) \in Gr(A) \iff \langle x_1-x_2, y_1 -y_2 \rangle \geq 0, \forall (x_2,y_2) \in Gr(A).
$$

We give an example to explain maximal monotone operators.

\bx\label{exmmo2}
For a lower semicontinuous convex function $\varphi:\mR^d\mapsto(-\infty, +\infty]$, we assume ${\rm Int}(Dom(\varphi))\neq \emptyset$, where $Dom(\varphi)\equiv\{x\in\mR^d; \varphi(x)<\infty\}$ and $\operatorname{Int}(Dom(\varphi))$ is the interior of $Dom(\varphi)$. Define the subdifferential operator of the function $\varphi$:
$$
\partial\varphi(x):=\{y\in\mR^d: \<y,z-x\>+\varphi(x)\leq \varphi(z), \forall z\in\mR^d\}.
$$
Then $\partial\varphi$ is a maximal monotone operator.
\ex

Take any $T>0$ and fix it. Let $\sV_{0}$ be the set of all continuous functions $K: [0,T]\mapsto\mR^d$ with finite variations and $K_{0} = 0$. For $K\in\sV_0$ and $s\in [0,T]$, we shall use $|K|_{0}^{s}$ to denote the variation of $K$ on [0,s]
and write $|K|_{TV}:=|K|_{0}^{T}$. Set
\ce
&&\sA:=\Big\{(X,K): X\in C([0,T],\overline{\cD(A)}), K \in \sV_0, \\
&&\qquad\qquad\quad~\mbox{and}~\langle X_{t}-x, \dif K_{t}-y\dif t\rangle \geq 0 ~\mbox{for any}~ (x,y)\in Gr(A)\Big\}.
\de
And about $\sA$ we have the following two results (cf.\cite{cepa2, ZXCH}).

\bl\label{equi}
For $X\in C([0,T],\overline{\cD(A)})$ and $K\in \sV_{0}$, the following statements are equivalent:
\begin{enumerate}[(i)]
	\item $(X,K)\in \sA$.
	\item For any $x, y\in C([0,T],\mR^d)$ with $(x_t, y_t)\in Gr(A)$, it holds that 
	$$
	\left\langle X_t-x_t, \dif K_t-y_t\dif t\right\rangle \geq0.
	$$
	\item For any $(X^{'},K^{'})\in \sA$, it holds that 
	$$
	\left\langle X_t-X_t^{'},\dif K_t-\dif K_t^{'}\right\rangle \geq0.
	$$
\end{enumerate}
\el

\bl\label{inteineq}
Assume that ${\rm Int}(\cD(A))\ne\emptyset$. For any $a\in {\rm Int}(\cD(A))$, there exist $M_1 >0$, and $M_{2},M_{3}\geq0$ such that  for any $(X,K)\in \sA$ and $0\leq s<t\leq T$,
$$
\int_s^t{\left< X_r-a, \dif K_r \right>}\geq M_1\left| K \right|_{s}^{t}-M_2\int_s^t{\left| X_r-a\right|}\dif r-M_3\left( t-s \right) .
$$
\el

\subsection{Backward stochastic variation inequalities}\label{bmsde}

In the subsection, we introduce backward SVIs. 

Consider the following backward SVI on $\mR^d$:
\be\left\{\begin{array}{l}
\dif Y_t \in \p\varphi(Y_t)\dif t -F(t,Y_t,Z_t)\dif t + Z_t\dif W_t,\\
Y_T=\xi, 
\end{array}
\right.
\label{bmsdeeq}
\ee
where $\varphi: \mR^d\rightarrow(-\infty, +\infty]$ is a proper ($\varphi\not\equiv+\infty$), convex and lower semicontinuous function, $\varphi(y)\geq \varphi(0)=0$, the coefficient $F: \Omega\times [0,T] \times \mR^d\times \mR^{d\times l}\mapsto{\mR^d}$ is Borel measurable, $\forall (y,z)\in\mR^d\times \mR^{d\times l}$, $F(\cdot,y,z)$ is $(\sF_t)_{t\in[0,T]}$-progressively measurable, and $\xi$ is a $\sF_T$-measurable random variable with values in $\overline{\cD\left( \p\varphi \right) }$ and $\mE[|\xi|^{2}+\varphi(\xi)]<\infty$. We define solutions for Eq.(\ref{bmsdeeq}).

\bd\label{solu}
We say that Eq.$(\ref{bmsdeeq})$ admits a solution with the terminal value $\xi$ if there exists a triple $\left\{ \left( Y_t, K_t, Z_t\right) :t\in \left[ 0,T \right] \right\} $ which is a $(\sF_t)_{t\in[0,T]}$-progressively measurable process and satisfies 
\begin{enumerate}[(i)]
\item 
$(Y_{\cdot},K_{\cdot})\in \sA$, $\dif \mP\times\dif t$-a.e. on $\Omega\times\left[ 0, T \right]$,
\item 
$$
\mathbb{E}\left( \underset{0\leqslant t\leqslant T}{\sup}\left| Y_t \right|^2+\int_0^T{\| Z_s \|^2}\dif s\right) <\infty, 
$$
\item 
\ce
Y_t=\xi -(K_T-K_t)+\int_t^T{g\left( s, Y_s, Z_s\right)}\dif s-\int_t^T{Z_s}\dif W_s.
\de
\end{enumerate}
\ed

\section{Main results}\label{ap}

In this section, we formulate the main results.

\subsection{The average principle for the SDEs}

In this subsection, we give an average principle for Eq.(\ref{multfsde}).

Consider Eq.(\ref{multfsde}), i.e.
\ce\left\{\begin{array}{l}
\dif X_{s}^{\e,t,\zeta}=b(\frac{s}{\e},X_{s}^{\e,t,\zeta})\dif s+\s(\frac{s}{\e},X_{s}^{\e,t,\zeta})\dif W_{s},\\
X_{t}^{\e,t,\zeta}=\zeta.
\end{array}
\right.
\de

We assume:
\begin{enumerate}[$(\mathbf{H}^1_{b, \s})$]
\item There exists a constant $L_1>0$ such that for any $s\in\mR_+, x_i\in\mR^m, i=1,2,$
\ce
&&|b(s,x_1)-b(s,x_2)|+\|\s(s,x_1)-\s(s,x_2)\|\leq L_1|x_1-x_2|,\\
&&|b(s,0)|+\|\s(s,0)\|\leq L_1.
\de
\end{enumerate}
\begin{enumerate}[$(\mathbf{H}^2_{b, \s})$]
\item 
There exist $\bar{b}: \mR^{m}\rightarrow\mR^{m}$, $\bar{\s}: \mR^{m}\rightarrow\mR^{m\times l}$ such that for any $\hat{T}\in\mR_+, x\in\mR^m$
\ce
&&\left|\frac{1}{\hat{T}}\int_0^{\hat{T}}b(s,x)\dif s-\bar{b}(x)\right|^2\leq \kappa_1(\hat{T})(1+|x|^2),\\
&&\frac{1}{\hat{T}}\int_0^{\hat{T}}\|\s(s,x)-\bar{\s}(x)\|^2\dif s\leq \kappa_2(\hat{T})(1+|x|^2),
\de
where $\kappa_i(\cdot)$ is a continuous and positive bounded function with $\lim\limits_{\hat{T}\rightarrow\infty}\kappa_i(\hat{T})=0, i=1,2$.
\end{enumerate}

\br\label{bs}
$\bar{b}, \bar{\s}$ are Lipschitz continuous in $x$. Indeed, for $x_1, x_2\in\mR^m$,
\ce
|\bar{b}(x_1)-\bar{b}(x_2)|^2&\leq& 3\left|\bar{b}(x_1)-\frac{1}{\hat{T}}\int_0^{\hat{T}}b(s,x_1)\dif s\right|^2+3\left|\frac{1}{\hat{T}}\int_0^{\hat{T}}b(s,x_1)\dif s-\frac{1}{\hat{T}}\int_0^{\hat{T}}b(s,x_2)\dif s\right|^2\\
&&+3\left|\frac{1}{\hat{T}}\int_0^{\hat{T}}b(s,x_2)\dif s-\bar{b}(x_2)\right|^2\\
&\leq& 3\left|\bar{b}(x_1)-\frac{1}{\hat{T}}\int_0^{\hat{T}}b(s,x_1)\dif s\right|^2+3\frac{1}{\hat{T}}\int_0^{\hat{T}}|b(s,x_1)-b(s,x_2)|^2\dif s\\
&&+3\left|\frac{1}{\hat{T}}\int_0^{\hat{T}}b(s,x_2)\dif s-\bar{b}(x_2)\right|^2\\
&\leq&C\kappa_1(\hat{T})(1+|x_1|^2+|x_2|^2)+3L^2_1|x_1-x_2|^2.
\de
Let $\hat{T}\rightarrow\infty$, and we obtain that $|\bar{b}(x_1)-\bar{b}(x_2)|^2\leq 3L^2_1|x_1-x_2|^2$. By the same deduction, one can verify that $\bar{\s}$ is Lipschitz continuous in $x$.
\er

Under $(\mathbf{H}^1_{b, \s})$, by \cite[Theorem 19.3]{h}, Eq.(\ref{multfsde}) has a unique solution $X^{\e,t,\zeta}$. Then we construct the following SDE:
\be\left\{\begin{array}{l}
\dif \bar{X}_{s}^{t,\zeta}=\bar{b}(\bar{X}_{s}^{t,\zeta})\dif s+\bar{\s}(\bar{X}_{s}^{t,\zeta})\dif W_{s},\\
\bar{X}_{t}^{t,\zeta}=\zeta.
\end{array}
\right.
\label{averfsde}
\ee
By Remark \ref{bs} and \cite[Theorem 19.3]{h}, we have that Eq.(\ref{averfsde}) also has a unique solution $\bar{X}^{t,\zeta}$. And the following theorem indicates the relationship between $X^{\e,t,\zeta}$ and $\bar{X}^{t,\zeta}$.

\bt\label{xbarxdiff}
Suppose that $(\mathbf{H}^1_{b, \s}), (\mathbf{H}^2_{b, \s})$ hold and $\mE|\zeta|^{2p+2}<\infty$ for any $p\geq 1$. Then it holds that
\be
\lim\limits_{\e\rightarrow 0}\mE\(\sup_{s\in[t,T]}|X_{s}^{\e,t,\zeta}-\bar{X}_{s}^{t,\zeta}|^{2p}\)=0.
\label{xebarx2p}
\ee
In particular, we have that for $0<\g<1$
\be
\mE\sup\limits_{s\in[t,T]}|X^{\e,t,\zeta}_s-\bar{X}^{t,\zeta}_{s}|^2\leq C(\e^\g+\e^{2\g}+\kappa_1(\e^{\g-1})+\kappa_2(\e^{\g-1})),
\label{xebarx2}
\ee
where the constant $C>0$ is independent of $\e$.
\et

The proof of the above theorem is placed in Section \ref{xbarxdiffproo}.

\br
In \cite{gxwh}, Guo et al. also studied Eq.(\ref{multfsde}) and obtained that 
$$
\lim\limits_{\e\rightarrow 0}\mE\(\sup_{s\in[t,T]}|X_{s}^{\e,t,\zeta}-\bar{X}_{s}^{t,\zeta}|^{2}\)=0
$$
under the local Lipschitz and monotone conditions. It is obvious that our result is better. 
\er

\subsection{The average principle for backward SVIs}

In this subsection, we present an average principle for Eq.(\ref{multfbsde}). 

Consider Eq.(\ref{multfbsde}), i.e.
\ce\left\{\begin{array}{l}
\dif Y_{s}^{\e,t,\zeta}\in \p\varphi(Y_{s}^{\e,t,\zeta})\dif s-[f_1(\frac{s}{\e},X_{s}^{\e,t,\zeta},Y_{s}^{\e,t,\zeta})+f_2(Z_{s}^{\e,t,\zeta})]\dif s+Z_{s}^{\e,t,\zeta}\dif W_{s},\\
Y_{T}^{\e,t,\zeta}=g(X_{T}^{\e,t,\zeta})\in\overline{\cD(\p\varphi)}.
\end{array}
\right.
\de

We assume:
\begin{enumerate}[$({\bf H}_{\varphi})$]
\item There exists a constant $L_2>0$ and a $q_1\in\mN\cup\{0\}$ 
such that 
$$
|\varphi(g(x))|\leq L_2(1+|x|^{q_1}), \quad x\in\mR^m.
$$
\end{enumerate}
\begin{enumerate}[(${\bf H}_{g}$)]
\item There exists a constant $L_3>0$ and a $q_2\in\mN$ such that for any $x_i\in\mR^m, i=1,2,$
\ce
|g(x_1)-g(x_2)|\leq L_3(1+|x_1|^{q_2}+|x_2|^{q_2})|x_{1}-x_{2}|.
\de
\end{enumerate}
\begin{enumerate}[$(\mathbf{H}^1_{f})$]
\item
There exist two constants $L_4>0, 0<L_5<1$ and a $q_3\in\mN$ such that for any $s\in\mR_+, x_i\in\mR^m, y_i\in\mR^d, z_i\in\mR^{d\times l}, i=1,2,$
\ce
&&|f_1(s,x_{1},y_{1})-f_1(s,x_{2},y_2)|\leq L_4((1+|x_1|^{q_3}+|x_2|^{q_3})|x_{1}-x_{2}|+|y_{1}-y_{2}|),\\
&&|f_1(s,0,0)|\leq L_4,\quad |f_2(z_1)-f_2(z_2)|\leq L_5\|z_1-z_2\|.
\de
\end{enumerate}

\begin{enumerate}[$(\mathbf{H}^2_{f})$]
\item
There exists $\bar{f}_1: \mR^{m}\times\mR^{d}\rightarrow\mR^{d}$ satisfying for any $\hat{T}\in\mR_+, x\in\mR^m, y\in\mR^d$
\ce
\left|\frac{1}{\hat{T}}\int_0^{\hat{T}}f_1(s,x,y)\dif s-\bar{f}_1(x,y)\right|^2\leq \kappa_3(\hat{T})(1+|x|^2+|y|^2),
\de
where $\kappa_3(\cdot)$ is a continuous and positive bounded function with $\lim\limits_{\hat{T}\rightarrow\infty}\kappa_3(\hat{T})=0$.
\end{enumerate}

\br\label{con}
$(i)$ By $({\bf H}_{g})$, it holds that for $x\in\mR^m$
\be
|g(x)|\leq (2L_3+|g(0)|)(1+|x|^{q_2+1}).
\label{glingro}
\ee

$(ii)$ $(\mathbf{H}^1_{f})$ implies that for $s\in\mR_+, x\in\mR^m, y\in\mR^d, z\in\mR^{d\times l}$
\be
|f_1(s,x,y)|\leq 2L_4(1+|x|^{q_3+1}+|y|), \quad |f_2(z)|\leq (L_5+|f_2(0)|)(1+\|z\|).
\label{flingro}
\ee

$(iii)$ By the similar deduction to that in Remark \ref{bs}, we know that $\bar{f}_1$ is Lipschitz continuous in $x, y$.
\er

Under $({\bf H}_{\varphi})$, $(\mathbf{H}_{g})$ and $(\mathbf{H}^1_{f})$, by \cite[Theorem 1.1]{pr}, the system (\ref{multfbsde}) has a unique solution $(Y^{\e,t,\zeta}, K^{\e,t,\zeta}, Z^{\e,t,\zeta})$. Then we construct the following backward SVI:
\be\left\{\begin{array}{l}
\dif \bar{Y}_{s}^{t,\zeta}\in \p\varphi(\bar{Y}_{s}^{t,\zeta})\dif s-[\bar{f}_1(\bar{X}_{s}^{t,\zeta},\bar{Y}_{s}^{t,\zeta})+f_2(\bar{Z}_{s}^{t,\zeta})]\dif s+\bar{Z}_{s}^{t,\zeta}\dif W_{s},\\
\bar{Y}_{T}^{t,\zeta}=g(\bar{X}_{T}^{t,\zeta})\in\overline{\cD(\p\varphi)}.
\end{array}
\right.
\label{averfbsde}
\ee

Now, it is the position to state the main result in this subsection.

\bt\label{averprin}
Assume that $(\mathbf{H}^{1}_{b, \s})$, $({\bf H}_{\varphi})$, $(\mathbf{H}_{g})$, $(\mathbf{H}^{1}_{f})$, $(\mathbf{H}^{2}_{b, \s})$, $(\mathbf{H}^{2}_{f})$ hold and $\mE|\zeta|^{2p+2}<\infty$ for any $p\geq 1$. Then it holds that
\ce
\lim\limits_{\e\rightarrow0}\mE\sup\limits_{s\in[t,T]}|Y^{\e,t,\zeta}_s-\bar{Y}^{t,\zeta}_{s}|^2=0,
\de
where $(\bar{Y}^{t,\zeta}, \bar{K}^{t,\zeta}, \bar{Z}^{t,\zeta})$ is a solution of Eq.(\ref{averfbsde}).
\et

The proof of the above theorem is placed in Section \ref{averprinproo}.

\br\label{f1sf20}
We mention that if $f_1(s,x,y)$ is independent of $s$ and $f_2(z)=0$, the backward equation in the system (\ref{multfbsde}) is the same to the equation (2.6) in \cite{eo}. There under similar assumptions Essaky and Ouknine proved that if $X^{\e,t,\zeta}$ converges in law to $\bar{X}^{t,\zeta}$, $Y^{\e,t,\zeta}$ also converges in law to $\bar{Y}^{t,\zeta}$ (cf. \cite[Theorem 3.1]{eo}). Since the mean square convergence implies the convergence in law, our result is stronger.
\er

\section{Proof of Theorem \ref{xbarxdiff}}\label{xbarxdiffproo}

In this section, we show Theorem \ref{xbarxdiff}. We begin with some key estimates.

\bl
Under $(\mathbf{H}^1_{b, \s})$ and $\mE|\zeta|^{2p+2}<\infty$ for any $p\geq 1$, it holds that
\be
&&\mE\sup\limits_{s\in[t,T]}|X^{\e,t,\zeta}_s|^{2p+2}\leq C(1+\mE|\zeta|^{2p+2}),\label{xesti}\\
&&\mE|X^{\e,t,\zeta}_{s+h}-X^{\e,t,\zeta}_s|^{2p+2}\leq C(h^{p+1}+h^{2p+2}), \quad t\leq s\leq s+h\leq T,\label{xeshs} 
\ee
where $C$ is independent of $\e$.
\el

\bl
Suppose that $(\mathbf{H}^1_{b, \s})$ and $(\mathbf{H}^2_{b, \s})$ hold and $\mE|\zeta|^{2p+2}<\infty$ for any $p\geq 1$. Then it holds that 
\be
\mE\sup\limits_{s\in[t,T]}|\bar{X}^{t,\zeta}_{s}|^{2p+2}\leq C(1+\mE|\zeta|^{2p+2}),
\label{barxesti}
\ee
where $C$ is independent of $\e$.
\el

Since the proofs of two above lemmas are standard, we omit them.

{\bf Proof of Theorem \ref{xbarxdiff}.}
First of all, by (\ref{multfsde}) and (\ref{averfsde}), it holds that
\be
&&\mE\sup\limits_{s\in[t,T]}|X^{\e,t,\zeta}_s-\bar{X}^{t,\zeta}_{s}|^2\no\\
&\leq& 2\mE\sup\limits_{s\in[t,T]}\left|\int_s^T\left(b(\frac{r}{\e},X^{\e,t,\zeta}_r)-\bar{b}(\bar{X}^{t,\zeta}_{r})\right)\dif r\right|^2+2C\mE\int_t^T\|\s(\frac{r}{\e},X^{\e,t,\zeta}_r)-\bar{\s}(\bar{X}^{t,\zeta}_{r})\|^2\dif r\no\\
&\leq& 4\mE\sup\limits_{s\in[t,T]}\left|\int_s^T\left(b(\frac{r}{\e},X^{\e,t,\zeta}_r)-\bar{b}(X^{\e,t,\zeta}_r)\right)\dif r\right|^2+4\mE\sup\limits_{s\in[t,T]}\left|\int_s^T\left(\bar{b}(X^{\e,t,\zeta}_r)-\bar{b}(\bar{X}^{t,\zeta}_{r})\right)\dif r\right|^2\no\\
&&+4C\mE\int_t^T\|\s(\frac{r}{\e},X^{\e,t,\zeta}_r)-\bar{\s}(X^{\e,t,\zeta}_r)\|^2\dif r+4C\mE\int_t^T\|\bar{\s}(X^{\e,t,\zeta}_r)-\bar{\s}(\bar{X}^{t,\zeta}_{r})\|^2\dif r\no\\
&\leq&C\int_t^T\mE\sup\limits_{s\in[t,r]}|X^{\e,t,\zeta}_s-\bar{X}^{t,\zeta}_{s}|^2\dif r+ 4\mE\sup\limits_{s\in[t,T]}\left|\int_s^T\left(b(\frac{r}{\e},X^{\e,t,\zeta}_r)-\bar{b}(X^{\e,t,\zeta}_r)\right)\dif r\right|^2\no\\
&&+4C\mE\int_t^T\|\s(\frac{r}{\e},X^{\e,t,\zeta}_r)-\bar{\s}(X^{\e,t,\zeta}_r)\|^2\dif r\no\\
&=:&C\int_t^T\mE\sup\limits_{s\in[t,r]}|X^{\e,t,\zeta}_s-\bar{X}^{t,\zeta}_{s}|^2\dif r+J_1+J_2.
\label{bszong}
\ee

For $J_1$, we define $X^{\e,t,\zeta}_r=\zeta, r\in[0,t]$ and by $(\mathbf{H}^{1}_{b, \s})$ $(\mathbf{H}^2_{b, \s})$ get that
\ce
J_1&\leq& 8\mE\sup\limits_{s\in[0,T]}\left|\int_0^s\left(b(\frac{r}{\e},X^{\e,t,\zeta}_r)-\bar{b}(X^{\e,t,\zeta}_r)-b(\frac{r}{\e},X^{\e,t,\zeta}_{r(\d)})+\bar{b}(X^{\e,t,\zeta}_{r(\d)})\right)\dif r\right|^2\no\\
&&+8\mE\sup\limits_{s\in[0,T]}\left|\int_0^s\left(b(\frac{r}{\e},X^{\e,t,\zeta}_{r(\d)})-\bar{b}(X^{\e,t,\zeta}_{r(\d)})\right)\dif r\right|^2\no\\
&\leq& C\int_0^T\mE|X^{\e,t,\zeta}_r-X^{\e,t,\zeta}_{r(\d)}|^2\dif r+8\mE\sup\limits_{s\in[0,T]}\left|\int_0^{[\frac{s}{\d}]\d}\left(b(\frac{r}{\e},X^{\e,t,\zeta}_{r(\d)})-\bar{b}(X^{\e,t,\zeta}_{r(\d)})\right)\dif r\right|^2\no\\
&&+8\mE\sup\limits_{s\in[0,T]}\left|\int_{[\frac{s}{\d}]\d}^s\left(b(\frac{r}{\e},X^{\e,t,\zeta}_{r(\d)})-\bar{b}(X^{\e,t,\zeta}_{r(\d)})\right)\dif r\right|^2\no\\
&=:&J_{11}+J_{12}+J_{13},
\de
where $\d$ is a fixed positive number depending on $\e$, $r(\d):=[\frac{r}{\d}]\d$ and $[\frac{r}{\d}]$ denotes the integer part of $\frac{r}{\d}$.

For $J_{11}$, (\ref{xeshs}) implies that 
\ce
J_{11}\leq CT(\d+\d^2).
\de

For $J_{12}$, from $(\mathbf{H}^2_{b, \s})$ it follows that
\ce
J_{12}&\leq&8[\frac{T}{\d}]\sum\limits_{k=0}^{[\frac{T}{\d}]-1}\mE\left|\int_{k\d}^{(k+1)\d}\left(b(\frac{r}{\e},X^{\e,t,\zeta}_{k\d})-\bar{b}(X^{\e,t,\zeta}_{k\d})\right)\dif r\right|^2\\
&\leq&8T\d\sum\limits_{k=0}^{[\frac{T}{\d}]-1}\mE\left|\frac{\e}{\d}\int_{k\d/\e}^{(k+1)\d/\e}\left(b(u,X^{\e,t,\zeta}_{k\d})-\bar{b}(X^{\e,t,\zeta}_{k\d})\right)\dif u\right|^2\\
&\leq&8T^2\kappa_1(\frac{\d}{\e})\sup\limits_{0\leq k\leq [\frac{T}{\d}]-1}\mE(1+|X^{\e,t,\zeta}_{k\d}|^2)\\
&\leq&C\kappa_1(\frac{\d}{\e}).
\de

For $J_{13}$, by the linear growth of $b, \bar{b}$ and (\ref{xesti})  it holds that
\ce
J_{13}&\leq& 8\d\mE\sup\limits_{s\in[0,T]}\int_{[\frac{s}{\d}]\d}^s\left|b(\frac{r}{\e},X^{\e,t,\zeta}_{r(\d)})-\bar{b}(X^{\e,t,\zeta}_{r(\d)})\right|^2\dif r\\
&\leq&C\d\mE\int_{0}^T(1+|X^{\e,t,\zeta}_{r(\d)}|^2)\dif r\\
&\leq& C\d.
\de

Combining the above deduction, we obtain that
\be
J_1\leq C(\d+\d^2+\kappa_1(\frac{\d}{\e})).
\label{j1}
\ee
Then the same computation to that for $J_1$ yields that
\be
J_2\leq C(\d+\d^2+\kappa_2(\frac{\d}{\e})).
\label{j2}
\ee

Finally, inserting (\ref{j1}) and (\ref{j2}) into (\ref{bszong}), by the Gronwall inequality we have that
\ce
\mE\sup\limits_{s\in[t,T]}|X^{\e,t,\zeta}_s-\bar{X}^{t,\zeta}_{s}|^2\leq C(\d+\d^2+\kappa_1(\frac{\d}{\e})+\kappa_2(\frac{\d}{\e})).
\de
By taking $\d=\e^\g$ for $0<\g<1$, we obtain (\ref{xebarx2}). 

Next, by the Chebyshev inequality and (\ref{xebarx2}), it holds that for any $\t>0$
$$ 
\mP\(\sup_{s\in[t,T]}|X_{s}^{\e,t,\zeta}-\bar{X}_{s}^{t,\zeta}|>\t\)\leq \frac{\mE\(\sup_{s\in[t,T]}|X_{s}^{\e,t,\zeta}-\bar{X}_{s}^{t,\zeta}|^{2}\)}{\t^2}\leq \frac{C}{\t^2}\(\d+\d^2+\kappa_1(\e^{\g-1})+\kappa_2(\e^{\g-1})\),
$$
which implies that
$$
\sup_{s\in[t,T]}|X_{s}^{\e,t,\zeta}-\bar{X}_{s}^{t,\zeta}|\overset{\mP}{\rightarrow} 0,
$$
as $\e$ tends to $0$. Besides, from (\ref{xesti}) and (\ref{barxesti}) it follows that 
$$
\sup\limits_{\e}\mE\sup_{s\in[t,T]}|X_{s}^{\e,t,\zeta}-\bar{X}_{s}^{t,\zeta}|^{2p+2}\leq C(1+\mE|\zeta|^{2p+2}).
$$
Therefore, by the Vitali convergence theorem one can obtain that
$$
\lim\limits_{\e\rightarrow 0}\mE\(\sup_{s\in[t,T]}|X_{s}^{\e,t,\zeta}-\bar{X}_{s}^{t,\zeta}|^{2p}\)=0,
$$
which completes the proof.

\section{Proofs of Theorem \ref{averprin}}\label{averprinproo} 

In this section, we prove Theorem \ref{averprin}. First of all, we prepare some key lemmas.

\bl\label{xeyeze}
Suppose that $(\mathbf{H}^{1}_{b, \s})$, $(\mathbf{H}_{\varphi})$, $(\mathbf{H}_{g})$ and $(\mathbf{H}^{1}_{f})$ hold and $\mE|\zeta|^{2p+2}<\infty$ for any $p\geq 1$. Then there exists a constant $C>0$ such that 
\be
\mE\sup\limits_{s\in[t,T]}|Y_s^{\e,t,\zeta}|^2+\int_t^T\mE\|Z_r^{\e,t,\zeta}\|^2\dif r+\mE|K^{\e,t,\zeta}|_t^T\leq C(1+\mE|\zeta|^{2q_2+2}+\mE|\zeta|^{2q_3+2}).
\label{yzetx}
\ee
\el
\begin{proof}
By the It\^o formula, it holds that 
\be
&&|Y_s^{\e,t,\zeta}|^2+\int_s^T\|Z_r^{\e,t,\zeta}\|^2\dif r\no\\
&=&|g(X_T^{\e,t,\zeta})|^2-2\int_s^T\<Y_r^{\e,t,\zeta},\dif K_r^{\e,t,\zeta}\>-2\int_s^T\<Y_r^{\e,t,\zeta},Z_r^{\e,t,\zeta}\dif W_r\>\no\\
&&+2\int_s^T\<Y_r^{\e,t,\zeta}, f_1(\frac{r}{\e}, X_r^{\e,t,\zeta}, Y_r^{\e,t,\zeta})+f_2(Z_r^{\e,t,\zeta})\>\dif r.
\label{yzito}
\ee
Taking the expectation on two sides, by Lemma \ref{equi} and $(\ref{glingro})$ $(\ref{flingro})$, we know that for any $v\in\p\varphi(0)$,
\ce
\mE|Y_s^{\e,t,\zeta}|^2+\int_s^T\mE\|Z_r^{\e,t,\zeta}\|^2\dif r&\leq&C(1+\mE|\zeta|^{2q_2+2})+2\int_s^T\mE|Y_r^{\e,t,\zeta}||v|\dif r+C\int_s^T\mE|Y_r^{\e,t,\zeta}|^2\dif r\\
&&+CT+C\int_s^T\mE|X_r^{\e,t,\zeta}|^{2q_3+2}\dif r+\frac{1}{2}\int_s^T\mE\|Z_r^{\e,t,\zeta}\|^2\dif r,
\de
and
\be
\mE|Y_s^{\e,t,\zeta}|^2+\frac{1}{2}\int_s^T\mE\|Z_r^{\e,t,\zeta}\|^2\dif r&\leq& C(1+\mE|\zeta|^{2q_2+2}+\mE|\zeta|^{2q_3+2})+(|v|^2+C)T\no\\
&&+C\int_s^T\mE|Y_r^{\e,t,\zeta}|^2\dif r.
\label{yeze}
\ee
The Gronwall inequality implies that
\ce
\sup\limits_{s\in[t,T]}\mE|Y_s^{\e,t,\zeta}|^2\leq C(1+\mE|\zeta|^{2q_2+2}+\mE|\zeta|^{2q_3+2}),
\de
which together with (\ref{yeze}) yields that
\ce
\int_t^T\mE\|Z_r^{\e,t,\zeta}\|^2\dif r\leq C(1+\mE|\zeta|^{2q_2+2}+\mE|\zeta|^{2q_3+2}).
\de

In the following, we investigate (\ref{yzito}) and by the BDG inequality obtain that
\ce
\mE\sup\limits_{s\in[t,T]}|Y_s^{\e,t,\zeta}|^2&\leq&C(1+\mE|\zeta|^{2q_2+2})+2\int_t^T\mE|Y_r^{\e,t,\zeta}||v|\dif r+C\mE\left(\int_t^T|Y_r^{\e,t,\zeta}|^2\|Z_r^{\e,t,\zeta}\|^2\dif r\right)^{1/2}\\
&&+CT+C\int_t^T\mE|X_r^{\e,t,\zeta}|^{2q_3+2}\dif r+C\int_t^T\mE|Y_r^{\e,t,\zeta}|^2\dif r+C\int_t^T\mE\|Z_r^{\e,t,\zeta}\|^2\dif r\\
&\leq&C(1+\mE|\zeta|^{2q_2+2}+\mE|\zeta|^{2q_3+2})+C\int_t^T\mE\sup\limits_{s\in[t,r]}|Y_s^{\e,t,\zeta}|^2\dif r+\frac{1}{2}\mE\sup\limits_{r\in[t,T]}|Y_r^{\e,t,\zeta}|^2,
\de
which together with the Gronwall inequality yields that
\ce
\mE\sup\limits_{r\in[t,T]}|Y_r^{\e,t,\zeta}|^2\leq C(1+\mE|\zeta|^{2q_2+2}+\mE|\zeta|^{2q_3+2}).
\de

Finally, based on Lemma \ref{inteineq} and (\ref{yzito}), it holds that
\ce
2M_1|K^{\e,t,\zeta}|_t^T&\leq&|g(X_T^{\e,t,\zeta})|^2+2M_2\int_t^T|Y_s^{\e,t,\zeta}|\dif s+2M_3T-2\int_t^T\<Y_s^{\e,t,\zeta},Z_s^{\e,t,\zeta}\dif W_s\>\\
&&+2\int_t^T\<Y_s^{\e,t,\zeta}, f_1(\frac{s}{\e}, X_s^{\e,t,\zeta}, Y_s^{\e,t,\zeta})+f_2(Z_s^{\e,t,\zeta})\>\dif s.
\de
Hence by $(\ref{glingro})$ $(\ref{flingro})$
\ce
\mE|K^{\e,t,\zeta}|_t^T\leq C(1+\mE|\zeta|^{2q_2+2}+\mE|\zeta|^{2q_3+2}).
\de
The proof is complete.
\end{proof}

\bl\label{xeyeze}
Suppose that $(\mathbf{H}^{1}_{b, \s})$, $(\mathbf{H}_{\varphi})$, $(\mathbf{H}_{g})$, $(\mathbf{H}^{1}_{f})$, $(\mathbf{H}^{2}_{b, \s})$, $(\mathbf{H}^{2}_{f})$ hold and $\mE|\zeta|^{2p+2}<\infty$ for any $p\geq 1$. Then Eq.(\ref{averfbsde}) has a unique solution $(\bar{Y}^{t,\zeta}, \bar{K}^{t,\zeta}, \bar{Z}^{t,\zeta})$. Moreover, there exists a constant $C>0$ such that 
\be
\mE\sup\limits_{s\in[t,T]}|\bar{Y}^{t,\zeta}_s|^2+\int_t^T\mE\|\bar{Z}^{t,\zeta}_r\|^2\dif r+\mE|\bar{K}^{t,\zeta}|_t^T\leq C(1+\mE|\zeta|^{2q_2+2}+\mE|\zeta|^{2q_3+2}).
\label{baryztx}
\ee
\el
\begin{proof}
By Remark \ref{con} and \cite[Theorem 1.1]{pr}, Eq.(\ref{averfbsde}) has a unique solution $(\bar{Y}^{t,\zeta}, \bar{K}^{t,\zeta}, \bar{Z}^{t,\zeta})$. Then by the same deduction to that of (\ref{yzetx}), we can show (\ref{baryztx}).
\end{proof}

\bl
Under the assumptions of Theorem \ref{averprin}, it holds that 
\be
&&\lim\limits_{\varrho\rightarrow0}\sup\limits_{s\in[t,T]}\mE\sup _{s \leqslant r \leqslant s+\varrho}|Y_{r}^{\e,t,\zeta}-Y_{s+\varrho}^{\e,t,\zeta}|^{2}=0, \label{yessro}\\
&&\lim\limits_{\varrho\rightarrow0}\sup\limits_{s\in[t,T]}\mE\sup _{s \leqslant r \leqslant s+\varrho}|\bar{Y}_{r}^{t,\zeta}-\bar{Y}_{s+\varrho}^{t,\zeta}|^{2}=0.\label{baryssro}
\ee
\el
\begin{proof}
Since the proofs of (\ref{yessro}) and (\ref{baryssro}) are similar, we only prove (\ref{yessro}). 

First of all, we know that for $\varrho>0$ and $t\leq s\leq r\leq s+\varrho\leq T$,
\ce
Y_r^{\e,t,\zeta}=Y_{s+\varrho}^{\e,t,\zeta}-K^{\e,t,\zeta}_{s+\varrho}+K^{\e,t,\zeta}_r+\int_r^{s+\varrho}[f_1(\frac{u}{\e},X_u^{\e,t,\zeta},Y_u^{\e,t,\zeta})+f_2(Z_u^{\e,t,\zeta})]\dif u-\int_r^{s+\varrho}Z_u^{\e,t,\zeta}\dif W_u.
\de
Then by the It\^o formula it holds that
\be
&&|Y_r^{\e,t,\zeta}-Y_{s+\varrho}^{\e,t,\zeta}|^2+\int_r^{s+\varrho}\|Z_u^{\e,t,\zeta}\|^2\dif u\no\\
&=&-2\int_r^{s+\varrho}\<Y_u^{\e,t,\zeta}-Y_{s+\varrho}^{\e,t,\zeta},\dif K^{\e,t,\zeta}_u\>+2\int_r^{s+\varrho}\<Y_u^{\e,t,\zeta}-Y_{s+\varrho}^{\e,t,\zeta},f_1(\frac{u}{\e},X_u^{\e,t,\zeta},Y_u^{\e,t,\zeta})\>\dif u\no\\
&&+2\int_r^{s+\varrho}\<Y_u^{\e,t,\zeta}-Y_{s+\varrho}^{\e,t,\zeta},f_2(Z_u^{\e,t,\zeta})\>\dif u-2\int_r^{s+\varrho}\<Y_u^{\e,t,\zeta}-Y_{s+\varrho}^{\e,t,\zeta},Z_u^{\e,t,\zeta}\dif W_u\>.
\label{yetsroito}
\ee

Next, we compute $-2\int_r^{s+\varrho}\<Y_u^{\e,t,\zeta}-Y_{s+\varrho}^{\e,t,\zeta},\dif K^{\e,t,\zeta}_u\>$. Take any $a\in{\rm Int}(\cD(\p\varphi))$. Thus, there is a $\t_0>0$ such that for any $R>0$ and $0<\t<\t_0$,
$$
\left\{x \in B(a,R): d(x,(\overline{\cD(\p\varphi)})^c) \geqslant \t\right\} \neq \emptyset,
$$
where $B(a,R):=\{x\in\mR^d: |x-a|\leq R\}$, $d(\cdot,\cdot)$ is the Euclidean distance in $\mR^d$ and $(\overline{\cD(\p\varphi)})^c$ denotes the complement of $\overline{\cD(\p\varphi)}$. Set
$$
g_R(\t):=\sup \left\{|z|: z \in \p\varphi(x) \text { for all } x \in B(a,R) \text { with } d\left(x,(\overline{\cD(\p\varphi)})^c\right) \geqslant \t\right\},
$$
and by the local boundedness of $\p\varphi$ on ${\rm Int}(\cD(\p\varphi))$, it holds that
$$
g_R(\t)<+\infty.
$$
Again put 
$$
h_R(\varrho):=\inf \left\{\t \in\left(0, \t_0\right): g_R(\t) \leqslant \varrho^{-1 / 2}\right\}, \quad \varrho>0,
$$
and we have that
$$
g_R\left(\varrho+h_R(\varrho)\right) \leqslant \varrho^{-1 / 2} \text { and } \quad \lim _{\varrho \downarrow 0} h_R(\varrho)=0.
$$
Take $\varrho_R>0$ such that $\varrho_R+h_R\left(\varrho_R\right)<\t_0$. For $0<\varrho<\varrho_R \wedge 1$, let $Y_{s+\varrho}^{\varepsilon, \varrho, R}$ be the projection of $Y_{s+\varrho}^{\e,t,\zeta}$ on $\left\{x \in B(a,R): d\left(x,(\overline{\cD(\p\varphi)})^c\right) \geqslant \varrho+h_R(\varrho)\right\}$. 
Thus, for $\Pi_{s+\varrho}^{\varepsilon, \varrho, R} \in \p\varphi(Y_{s+\varrho}^{\varepsilon,\varrho, R}), \sup\limits_{v\in[t,T]}|Y_v^{\e,t,\zeta}-a| \leqslant R$ and $0<r-s<\varrho$, it holds that
\ce
&&-2\int_r^{s+\varrho}\<Y_u^{\e,t,\zeta}-Y_{s+\varrho}^{\e,t,\zeta},\dif K^{\e,t,\zeta}_u\>\\
 &=&-2 \int_r^{s+\varrho}\left\langle Y_u^{\e,t,\zeta}-Y_{s+\varrho}^{\varepsilon, \varrho, R}, \dif K_u^{\e,t,\zeta} \right\rangle-2 \int_r^{s+\varrho}\left\langle Y_{s+\varrho}^{\varepsilon, \varrho, R}-Y_{s+\varrho}^{\e,t,\zeta}, \dif K_u^{\e,t,\zeta} \right\rangle \\ 
&\leqslant& -2 \int_r^{s+\varrho}\left\langle Y_u^{\e,t,\zeta}-Y_{s+\varrho}^{\varepsilon, \varrho, R}, \Pi_{s+\varrho}^{\varepsilon, \varrho, R}\right\rangle \dif u+2\left(\varrho+h_R(\varrho)\right)\left|K^{\e,t,\zeta} \right|_t^T\\
&\leqslant& 4 \varrho^{1 / 2}(R+|a|)+2\left(\varrho+h_R(\varrho)\right)\left|K^{\e,t,\zeta} \right|_t^T,
\de
and furthermore by (\ref{yetsroito})
\ce
&&\sup _{s \leqslant r \leqslant s+\varrho}\mE|Y_r^{\e,t,\zeta}-Y_{s+\varrho}^{\e,t,\zeta}|^2I_{\{\sup\limits_{v\in[t,T]}|Y_v^{\e,t,\zeta}-a| \leqslant R\}}+\mE\int_s^{s+\varrho}\|Z_u^{\e,t,\zeta}\|^2\dif uI_{\sup\limits_{v\in[t,T]}|Y_v^{\e,t,\zeta}-a| \leqslant R}\\
&\leq&4 \varrho^{1 / 2}(R+|a|)+2\left(\varrho+h_R(\varrho)\right)\mE\left|K^{\e,t,\zeta} \right|_t^T+C\mE\int_s^{s+\varrho}|Y_u^{\e,t,\zeta}-Y_{s+\varrho}^{\e,t,\zeta}|^2I_{\{\sup\limits_{v\in[t,T]}|Y_v^{\e,t,\zeta}-a| \leqslant R\}}\dif u\\
&&+C\varrho+C\int_s^{s+\varrho}\mE|X_u^{\e,t,\zeta}|^{2q_3+2}\dif u+C\int_s^{s+\varrho}\mE|Y_u^{\e,t,\zeta}|^2\dif u+\frac{1}{2}\mE\int_s^{s+\varrho}\|Z_u^{\e,t,\zeta}\|^2\dif uI_{\sup\limits_{v\in[t,T]}|Y_v^{\varepsilon}-a| \leqslant R}.
\de
Thus by the Gronwall inequality we obtain that
\ce
&&\sup _{s \leqslant r \leqslant s+\varrho}\mE|Y_r^{\e,t,\zeta}-Y_{s+\varrho}^{\e,t,\zeta}|^2I_{\{\sup\limits_{v\in[t,T]}|Y_v^{\e,t,\zeta}-a| \leqslant R\}}\\
&\leq& \(4 \varrho^{1 / 2}(R+|a|)+2\left(\varrho+h_R(\varrho)\right)\mE\left|K^{\e,t,\zeta} \right|_t^T\)+C\varrho,\\
&&\mE\int_s^{s+\varrho}\|Z_u^{\e,t,\zeta}\|^2\dif uI_{\sup\limits_{v\in[t,T]}|Y_v^{\varepsilon}-a| \leqslant R}\\
&\leq& 2\(4 \varrho^{1 / 2}(R+|a|)+2\left(\varrho+h_R(\varrho)\right)\mE\left|K^{\e,t,\zeta} \right|_t^T\)+C\varrho.
\de

Next, for (\ref{yetsroito}), from the BDG inequality and the H\"older inequality, it follows that
\ce
&&\mE\sup _{s \leqslant r \leqslant s+\varrho}|Y_r^{\e,t,\zeta}-Y_{s+\varrho}^\e|^2I_{\sup\limits_{v\in[t,T]}|Y_v^{\varepsilon}-a| \leqslant R}+\mE\int_s^{s+\varrho}\|Z_u^{\e,t,\zeta}\|^2\dif uI_{\sup\limits_{v\in[t,T]}|Y_v^{\varepsilon}-a| \leqslant R}\\
&\leq&4 \varrho^{1 / 2}(R+|a|)+2\left(\varrho+h_R(\varrho)\right)\mE\left|K^{\e,t,\zeta} \right|_t^T+C\int_s^{s+\varrho}\mE\sup _{u\leqslant\nu\leqslant s+\varrho}|Y_{\nu}^{\e,t,\zeta}-Y_{s+\varrho}^{\e,t,\zeta}|^2I_{\sup\limits_{v\in[t,T]}|Y_v^{\varepsilon}-a| \leqslant R}\dif u\\
&&+C\varrho+C\int_s^{s+\varrho}\mE|X_u^{\e,t,\zeta}|^{2q_3+2}\dif u+C\int_s^{s+\varrho}\mE|Y_u^{\e,t,\zeta}|^2\dif u+C\mE\int_s^{s+\varrho}\|Z_u^{\e,t,\zeta}\|^2\dif uI_{\sup\limits_{v\in[t,T]}|Y_v^{\varepsilon}-a| \leqslant R}\\
&&+\frac{1}{2}\mE\sup _{s \leqslant r \leqslant s+\varrho}|Y_r^{\e,t,\zeta}-Y_{s+\varrho}^{\e,t,\zeta}|^2I_{\sup\limits_{v\in[t,T]}|Y_v^{\varepsilon}-a| \leqslant R},
\de
which yields that
\ce
\mE\sup _{s \leqslant r \leqslant s+\varrho}|Y_r^{\e,t,\zeta}-Y_{s+\varrho}^{\e,t,\zeta}|^2I_{\sup\limits_{v\in[t,T]}|Y_v^{\varepsilon}-a| \leqslant R}\leq C\(4 \varrho^{1 / 2}(R+|a|)+2\left(\varrho+h_R(\varrho)\right)\mE\left|K^{\e,t,\zeta} \right|_t^T\)+C\varrho
\de
Based on this, letting $\varrho\rightarrow 0$ and $R\rightarrow \infty$, one could conclude (\ref{yessro}), which completes the proof.
\end{proof}

{\bf Proof of Theorem \ref{averprin}.}
{\bf Step 1.} We estimate $\int_t^T\mE\|Z_r^{\e,t,\zeta}-\bar{Z}_r^{t,\zeta}\|^2\dif r$.

First of all, combining (\ref{multfbsde}) with (\ref{averfbsde}), we have that
\ce
Y_s^{\e,t,\zeta}-\bar{Y}_s^{t,\zeta}&=&g(X_T^{\e,t,\zeta})-g(\bar{X}_T^{t,\zeta})-(K^{\e,t,\zeta}_T-K^{\e,t,\zeta}_s)+(\bar{K}^{t,\zeta}_T-\bar{K}^{t,\zeta}_s)\\
&&+\int_s^T\(f_1(\frac{r}{\e},X^{\e,t,\zeta}_r,Y_r^{\e,t,\zeta})-\bar{f}_1(\bar{X}^{t,\zeta}_r,\bar{Y}^{t,\zeta}_r)\)\dif r\\
&&+\int_s^T(f_2(Z_r^{\e,t,\zeta})-f_2(\bar{Z}_r^{t,\zeta}))\dif r-\int_s^T\(Z_r^{\e,t,\zeta}-\bar{Z}_r^{t,\zeta}\)\dif W_r.
\de
The It\^o formula for $|Y_s^{\e,t,\zeta}-\bar{Y}_s^{t,\zeta}|^2$ and $(\mathbf{H}_{g})$ $(\mathbf{H}^{1}_{f})$ imply that 
\be
&&|Y_s^{\e,t,\zeta}-\bar{Y}_s^{t,\zeta}|^2+\int_s^T\|Z_r^{\e,t,\zeta}-\bar{Z}_r^{t,\zeta}\|^2\dif r\no\\
&=&|g(X_T^{\e,t,\zeta})-g(\bar{X}_T^{t,\zeta})|^2-2\int_s^T\<Y_r^{\e,t,\zeta}-\bar{Y}_r^{t,\zeta},\dif (K^{\e,t,\zeta}_r-\bar{K}^{t,\zeta}_r)\>\no\\
&&-2\int_s^T\<Y_r^{\e,t,\zeta}-\bar{Y}_r^{t,\zeta},\(Z_r^{\e,t,\zeta}-\bar{Z}_r^{t,\zeta}\)\dif W_r\>\no\\
&&+2\int_s^T\<Y_r^{\e,t,\zeta}-\bar{Y}_r^{t,\zeta},f_1(\frac{r}{\e},X^{\e,t,\zeta}_r,Y_r^{\e,t,\zeta})-\bar{f}_1(\bar{X}^{t,\zeta}_r,\bar{Y}^{t,\zeta}_r)\>\dif r\no\\
&&+2\int_s^T\<Y_r^{\e,t,\zeta}-\bar{Y}_r^{t,\zeta},f_2(Z_r^{\e,t,\zeta})-f_2(\bar{Z}_r^{t,\zeta})\>\dif r\no\\
&\overset{Lemma ~\ref{equi}}{\leq}&L^2_3(1+|X_T^{\e,t,\zeta}|^{q_2}+|\bar{X}_T^{t,\zeta}|^{q_2})^2|X_T^{\e,t,\zeta}-\bar{X}_T^{t,\zeta}|^2-2\int_s^T\<Y_r^{\e,t,\zeta}-\bar{Y}_r^{t,\zeta},\(Z_r^{\e,t,\zeta}-\bar{Z}_r^{t,\zeta}\)\dif W_r\>\no\\
&&+2\int_s^T\<Y_r^{\e,t,\zeta}-\bar{Y}_r^{t,\zeta},f_1(\frac{r}{\e},X^{\e,t,\zeta}_r,Y_r^{\e,t,\zeta})-\bar{f}_1(X^{\e,t,\zeta}_r,Y_r^{\e,t,\zeta})\>\dif r\no\\
&&+2\int_s^T\<Y_r^{\e,t,\zeta}-\bar{Y}_r^{t,\zeta},\bar{f}_1(X^{\e,t,\zeta}_r,Y_r^{\e,t,\zeta})-\bar{f}_1(\bar{X}^{t,\zeta}_r,\bar{Y}^{t,\zeta}_r)\>\dif r\no\\
&&+2\int_s^T\<Y_r^{\e,t,\zeta}-\bar{Y}_r^{t,\zeta},f_2(Z_r^{\e,t,\zeta})-f_2(\bar{Z}_r^{t,\zeta})\>\dif r\no\\
&\leq&L^2_3(1+|X_T^{\e,t,\zeta}|^{q_2}+|\bar{X}_T^{t,\zeta}|^{q_2})^2|X_T^{\e,t,\zeta}-\bar{X}_T^{t,\zeta}|^2-2\int_s^T\<Y_r^{\e,t,\zeta}-\bar{Y}_r^{t,\zeta},\(Z_r^{\e,t,\zeta}-\bar{Z}_r^{t,\zeta}\)\dif W_r\>\no\\
&&+2\int_s^T\<Y_r^{\e,t,\zeta}-\bar{Y}_r^{t,\zeta},f_1(\frac{r}{\e},X^{\e,t,\zeta}_r,Y_r^{\e,t,\zeta})-\bar{f}_1(X^{\e,t,\zeta}_r,Y_r^{\e,t,\zeta})\>\dif r\no\\
&&+6L^2_4\int_s^T(1+|X_r^{\e,t,\zeta}|^{q_3}+|\bar{X}_r^{t,\zeta}|^{q_3})^2|X_r^{\e,t,\zeta}-\bar{X}_r^{t,\zeta}|^2\dif r\no\\
&&+(2+6L^2_4)\int_s^T|Y_r^{\e,t,\zeta}-\bar{Y}_r^{t,\zeta}|^2\dif r+L_5\int_s^T\|Z_r^{\e,t,\zeta}-\bar{Z}_r^{t,\zeta}\|^2\dif r.
\label{yehatyito}
\ee
By taking the expectation on two sides, we have that
\ce
&&\mE|Y_s^{\e,t,\zeta}-\bar{Y}_s^{t,\zeta}|^2+\int_s^T\mE\|Z_r^{\e,t,\zeta}-\bar{Z}_r^{t,\zeta}\|^2\dif r\no\\
&\leq&C(1+\mE|\zeta|^{4q_2}+\mE|\zeta|^{4q_3})^{1/2}(\mE\sup\limits_{s\in[t,T]}|X_s^{\e,t,\zeta}-\bar{X}_s^{t,\zeta}|^4)^{1/2}\no\\
&&+(2+6L^2_4)\int_s^T\mE|Y_r^{\e,t,\zeta}-\bar{Y}_r^{t,\zeta}|^2\dif s+L_5\int_s^T\mE\|Z_r^{\e,t,\zeta}-\bar{Z}_r^{t,\zeta}\|^2\dif r\no\\
&&+2\mE\int_s^T\<Y_r^{\e,t,\zeta}-\bar{Y}_r^{t,\zeta}, f_1(\frac{r}{\e},X^{\e,t,\zeta}_r,Y_r^{\e,t,\zeta})-\bar{f}_1(X^{\e,t,\zeta}_r,Y_r^{\e,t,\zeta})\>\dif r.
\de
Note that $0<L_5<1$ and 
\be
&&2\sup\limits_{s\in[t,T]}\mE\left|\int_s^T\<Y_r^{\e,t,\zeta}-\bar{Y}_r^{t,\zeta},f_1(\frac{r}{\e},X^{\e,t,\zeta}_r,Y_r^{\e,t,\zeta})-\bar{f}_1(X^{\e,t,\zeta}_r,Y_r^{\e,t,\zeta})\>\dif r\right|\no\\
&\leq&C\left(\sup\limits_{s\in[t,T]}\mE\sup _{s \leq r \leq s+\d}|Y_r^{\e,t,\zeta}-Y_{s+\d}^{\e,t,\zeta}|^{2}+\sup\limits_{s\in[t,T]}\mE\sup _{s \leq r \leq s+\d}|\bar{Y}^{t,\zeta}_{r}-\bar{Y}^{t,\zeta}_{s+\d}|^{2}\right)^{1/2}\no\\
&&+C\sup\limits_{s\in[t,T]}\mE\sup _{s \leq r \leq s+\d}|Y_r^{\e,t,\zeta}-Y_{s+\d}^{\e,t,\zeta}|^{2}+4CT\kappa_3^{1/2}(\frac{\d}{\e})+C(\d+\d^2)\no\\
&&+2\int_t^T\sup\limits_{s\in[t,r]}\mE|Y_s^{\e,t,\zeta}-\bar{Y}_s^{t,\zeta}|^2\dif r,
\label{yyfbarf}
\ee
where $\d$ is the same to that in the proof of Theorem \ref{xbarxdiff}. Thus, it holds that
\be
&&\sup\limits_{s\in[t,T]}\mE|Y_s^{\e,t,\zeta}-\bar{Y}_s^{t,\zeta}|^2+(1-L_5)\int_t^T\mE\|Z_r^{\e,t,\zeta}-\bar{Z}_r^{t,\zeta}\|^2\dif r\no\\
&\leq&\Gamma(\e)+(4+6L^2_4)\int_t^T\sup\limits_{s\in[t,r]}\mE|Y_s^{\e,t,\zeta}-\bar{Y}_s^{t,\zeta}|^2\dif r,
\label{yzdiffesti}
\ee
where 
\ce
\Gamma(\e)&:=&C(1+\mE|\zeta|^{4q_2}+\mE|\zeta|^{4q_3})^{1/2}(\mE\sup\limits_{s\in[t,T]}|X_s^{\e,t,\zeta}-\bar{X}_s^{t,\zeta}|^4)^{1/2}\\
&&+C\sup\limits_{s\in[t,T]}\mE\sup _{s \leq r \leq s+\d}|Y_r^{\e,t,\zeta}-Y_{s+\d}^{\e,t,\zeta}|^{2}+4CT\kappa_3^{1/2}(\frac{\d}{\e})+C(\d+\d^2)\\
&&+C\left(\sup\limits_{s\in[t,T]}\mE\sup _{s \leq r \leq s+\d}|Y_r^{\e,t,\zeta}-Y_{s+\d}^{\e,t,\zeta}|^{2}+\sup\limits_{s\in[t,T]}\mE\sup _{s \leq r \leq s+\d}|\bar{Y}^{t,\zeta}_{r}-\bar{Y}^{t,\zeta}_{s+\d}|^{2}\right)^{1/2}.
\de
The Gronwall inequality implies that 
\ce
\sup\limits_{s\in[t,T]}\mE|Y_s^{\e,t,\zeta}-\bar{Y}_s^{t,\zeta}|^2\leq\Gamma(\e) e^{(4+6L^2_4)T}.
\de
Inserting the above inequality into (\ref{yzdiffesti}), we obtain that
\be
\int_t^T\mE\|Z_r^{\e,t,\zeta}-\bar{Z}_r^{t,\zeta}\|^2\dif r\leq C\Gamma(\e).
\label{zebarzes}
\ee

{\bf Step 2.} We prove $\lim\limits_{\e\rightarrow0}\mE\sup\limits_{s\in[t,T]}|Y^{\e,t,\zeta}_s-\bar{Y}^{t,\zeta}_{s}|^2=0$.

For (\ref{yehatyito}), by the BDG inequality it holds that
\ce
&&\mE\sup\limits_{s\in[t,T]}|Y_s^{\e,t,\zeta}-\bar{Y}_s^{t,\zeta}|^2+\int_t^T\mE\|Z_r^{\e,t,\zeta}-\bar{Z}_r^{t,\zeta}\|^2\dif r\\
&\leq&L^2_3\mE(1+|X_T^{\e,t,\zeta}|^{q_2}+|\bar{X}_T^{t,\zeta}|^{q_2})^2|X_T^{\e,t,\zeta}-\bar{X}_T^{t,\zeta}|^2\\
&&+C\mE\left(\int_t^T|Y_r^{\e,t,\zeta}-\bar{Y}_r^{t,\zeta}|^2\|Z_r^{\e,t,\zeta}-\bar{Z}_r^{t,\zeta}\|^2\dif r\right)^{1/2}\\
&&+2\mE\sup\limits_{s\in[t,T]}\left|\int_s^T\<Y_r^{\e,t,\zeta}-\bar{Y}_r^{t,\zeta},f_1(\frac{r}{\e},X^{\e,t,\zeta}_r,Y_r^{\e,t,\zeta})-\bar{f}_1(X^{\e,t,\zeta}_r,Y_r^{\e,t,\zeta})\>\dif r\right|\no\\
&&+6L^2_4\mE\int_t^T(1+|X_r^{\e,t,\zeta}|^{q_3}+|\bar{X}_r^{t,\zeta}|^{q_3})^2|X_r^{\e,t,\zeta}-\bar{X}_r^{t,\zeta}|^2\dif r\no\\
&&+(2+6L^2_4)\mE\int_t^T|Y_r^{\e,t,\zeta}-\bar{Y}_r^{t,\zeta}|^2\dif r+L_5\mE\int_t^T\|Z_r^{\e,t,\zeta}-\bar{Z}_r^{t,\zeta}\|^2\dif r\\
&\leq&C(1+\mE|\zeta|^{4q_2}+\mE|\zeta|^{4q_3})^{1/2}(\mE\sup\limits_{s\in[t,T]}|X_s^{\e,t,\zeta}-\bar{X}_s^{t,\zeta}|^4)^{1/2}+\frac{1}{2}\mE\sup\limits_{s\in[t,T]}|Y_s^{\e,t,\zeta}-\bar{Y}_s^{t,\zeta}|^2\\
&&+(2+6L^2_4)\int_t^T\mE\sup\limits_{s\in[t,r]}|Y_s^{\e,t,\zeta}-\bar{Y}_s^{t,\zeta}|^2\dif r+C\int_t^T\mE\|Z_r^{\e,t,\zeta}-\bar{Z}_r^{t,\zeta}\|^2\dif r\no\\
&&+2\mE\sup\limits_{s\in[t,T]}\left|\int_s^T\<Y_r^{\e,t,\zeta}-\bar{Y}_r^{t,\zeta},f_1(\frac{r}{\e},X^{\e,t,\zeta}_r,Y_r^{\e,t,\zeta})-\bar{f}_1(X^{\e,t,\zeta}_r,Y_r^{\e,t,\zeta})\>\dif r\right|.
\de

Besides, by the similar deduction to that for (\ref{yyfbarf}), we can obtain that
\ce
&&2\mE\sup\limits_{s\in[t,T]}\left|\int_s^T\<Y_r^{\e,t,\zeta}-\bar{Y}_r^{t,\zeta},f_1(\frac{r}{\e},X^{\e,t,\zeta}_r,Y_r^{\e,t,\zeta})-\bar{f}_1(X^{\e,t,\zeta}_r,Y_r^{\e,t,\zeta})\>\dif r\right|\\
&\leq&C\left(\sup\limits_{s\in[t,T]}\mE\sup _{s \leq r \leq s+\d}|Y_r^{\e,t,\zeta}-Y_{s+\d}^{\e,t,\zeta}|^{2}+\sup\limits_{s\in[t,T]}\mE\sup _{s \leq r \leq s+\d}|\bar{Y}^{t,\zeta}_{r}-\bar{Y}^{t,\zeta}_{s+\d}|^{2}\right)^{1/2}\no\\
&&+C\sup\limits_{s\in[t,T]}\mE\sup _{s \leq r \leq s+\d}|Y_r^{\e,t,\zeta}-Y_{s+\d}^{\e,t,\zeta}|^{2}+4CT\kappa_3^{1/2}(\frac{\d}{\e})+C(\d+\d^2)\no\\
&&+2\int_t^T\mE\sup\limits_{s\in[t,r]}|Y_s^{\e,t,\zeta}-\bar{Y}_s^{t,\zeta}|^2\dif r,
\de
which together with (\ref{zebarzes}) and the Gronwall inequality yields that
\ce
\mE\sup\limits_{s\in[t,T]}|Y_s^{\e,t,\zeta}-\bar{Y}_s^{t,\zeta}|^2\leq C\Gamma(\e).
\de

Finally, we study the limit of $\mE\sup\limits_{s\in[t,T]}|Y_s^{\e,t,\zeta}-\bar{Y}_s^{t,\zeta}|^2$ as $\e\rightarrow0$. On one hand, by (\ref{xebarx2p}) it holds that
\ce
\lim\limits_{\e\rightarrow0}\mE\sup\limits_{s\in[t,T]}|X^{\e,t,\zeta}_s-\bar{X}^{t,\zeta}_{s}|^4=0.
\de
On the other hand, we take $\d=\e^\g$ for $0<\g<1$ and by (\ref{yessro}), (\ref{baryssro}) have that
\ce
&&\lim\limits_{\e\rightarrow0}\sup\limits_{s\in[t,T]}\mE\sup _{s \leqslant r \leqslant s+\d}|Y_{r}^{\e,t,\zeta}-Y_{s+\d}^{\e,t,\zeta}|^{2}=0, \\
&&\lim\limits_{\e\rightarrow0}\sup\limits_{s\in[t,T]}\mE\sup _{s \leqslant r \leqslant s+\d}|\bar{Y}_{r}^{t,\zeta}-\bar{Y}_{s+\d}^{t,\zeta}|^{2}=0,\\
&&\lim\limits_{\e\rightarrow0}\kappa_3^{1/2}(\frac{\d}{\e})=0.
\de
Combining the above deduction, one can get that
\ce
\lim\limits_{\e\rightarrow0}\mE\sup\limits_{s\in[t,T]}|Y^{\e,t,\zeta}_s-\bar{Y}^{t,\zeta}_{s}|^2=0.
\de

{\bf Step 3.} We prove (\ref{yyfbarf}).

Set
$$
I:=2\sup\limits_{s\in[t,T]}\mE\left|\int_s^T\<Y_r^{\e,t,\zeta}-\bar{Y}_r^{t,\zeta},f_1(\frac{r}{\e},X^{\e,t,\zeta}_r,Y_r^{\e,t,\zeta})-\bar{f}_1(X^{\e,t,\zeta}_r,Y_r^{\e,t,\zeta})\>\dif r\right|,
$$
and 
\be
I&\leq&2\sup\limits_{s\in[t,T]}\mE\left|\int_s^T\<Y_r^{\e,t,\zeta}-\bar{Y}_r^{t,\zeta}-Y_{r(\d)}^{\e,t,\zeta}+\bar{Y}_{r(\d)}^{t,\zeta},f_1(\frac{r}{\e},X^{\e,t,\zeta}_r,Y_r^{\e,t,\zeta})-\bar{f}_1(X^{\e,t,\zeta}_r,Y_r^{\e,t,\zeta})\>\dif r\right|\no\\
&&+2\sup\limits_{s\in[t,T]}\mE\left|\int_s^T\<Y_{r(\d)}^{\e,t,\zeta}-\bar{Y}_{r(\d)}^{t,\zeta},f_1(\frac{r}{\e},X^{\e,t,\zeta}_r,Y_r^{\e,t,\zeta})-f_1(\frac{r}{\e},X^{\e,t,\zeta}_{r(\d)},Y_{r(\d)}^{\e,t,\zeta})\>\dif r\right|\no\\ 
&&+2\sup\limits_{s\in[t,T]}\mE\left|\int_s^T\<Y_{r(\d)}^{\e,t,\zeta}-\bar{Y}_{r(\d)}^{t,\zeta},\bar{f}_1(X^{\e,t,\zeta}_{r(\d)},Y^{\e,t,\zeta}_{r(\d)})-\bar{f}_1(X^{\e,t,\zeta}_r,Y_r^{\e,t,\zeta})\>\dif r\right|\no\\
&&+2\sup\limits_{s\in[t,T]}\mE\left|\int_s^T\<Y_{r(\d)}^{\e,t,\zeta}-\bar{Y}_{r(\d)}^{t,\zeta},f_1(\frac{r}{\e},X^{\e,t,\zeta}_{r(\d)},Y_{r(\d)}^{\e,t,\zeta})-\bar{f}_1(X^{\e,t,\zeta}_{r(\d)},Y^{\e,t,\zeta}_{r(\d)})\>\dif r\right|\no\\
&=:&I_1+I_2+I_3+I_4,
\label{i1234}
\ee

For $I_1$, by the H\"older inequality and the growth of $f_1, \bar{f}_1$, it holds that
\be
I_1&\leq& 2\left(\mE\int_t^T|Y_r^{\e,t,\zeta}-\bar{Y}_r^{t,\zeta}-Y_{r(\d)}^{\e,t,\zeta}+\bar{Y}_{r(\d)}^{t,\zeta}|^2\dif r\right)^{1/2}\no\\
&&\times\left(\mE\int_t^T|f_1(\frac{r}{\e},X^{\e,t,\zeta}_r,Y_r^{\e,t,\zeta})-\bar{f}_1(X^{\e,t,\zeta}_r,Y_r^{\e,t,\zeta})|^2\dif r\right)^{1/2}\no\\
&\leq& C\left(\int_t^T(\mE|Y_r^{\e,t,\zeta}-Y_{r(\d)}^{\e,t,\zeta}|^2+\mE|\bar{Y}_r^{t,\zeta}-\bar{Y}_{r(\d)}^{t,\zeta}|^2)\dif r\right)^{1/2}\no\\
&&\times\left(\int_t^T(1+\mE|X^{\e,t,\zeta}_r|^{2q_3+2}+\mE|Y_r^{\e,t,\zeta}|^2)\dif r\right)^{1/2}\no\\
&\leq& C\left(\sup\limits_{s\in[t,T]}\mE\sup _{s \leq r \leq s+\d}|Y_r^{\e,t,\zeta}-Y_{s+\d}^{\e,t,\zeta}|^{2}+\sup\limits_{s\in[t,T]}\mE\sup _{s \leq r \leq s+\d}|\bar{Y}^{t,\zeta}_{r}-\bar{Y}^{t,\zeta}_{s+\d}|^{2}\right)^{1/2}.
\label{i1}
\ee

For $I_2$, the H\"older inequality and $(\mathbf{H}^{1}_{f})$  imply that
\be
I_2&\leq&\int_t^T\mE|Y_{r(\d)}^{\e,t,\zeta}-\bar{Y}_{r(\d)}^{t,\zeta}|^2\dif r+\int_t^T\mE|f_1(\frac{r}{\e},X^{\e,t,\zeta}_r,Y_r^{\e,t,\zeta})-f_1(\frac{r}{\e},X^{\e,t,\zeta}_{r(\d)},Y_{r(\d)}^{\e,t,\zeta})|^2\dif r\no\\
&\leq&\int_t^T\sup\limits_{s\in[t,r]}\mE|Y_s^{\e,t,\zeta}-\bar{Y}_s^{t,\zeta}|^2\dif r+2L^2_4\int_t^T\mE|Y_r^{\e,t,\zeta}-Y_{r(\d)}^{\e,t,\zeta}|^2\dif r\no\\
&&+2L^2_4\int_t^T\mE(1+|X^{\e,t,\zeta}_r|^{q_3}+|X^{\e,t,\zeta}_{r(\d)}|^{q_3})^2|X^{\e,t,\zeta}_r-X^{\e,t,\zeta}_{r(\d)}|^2\dif r\no\\
&\leq&\int_t^T\sup\limits_{s\in[t,r]}\mE|Y_s^{\e,t,\zeta}-\bar{Y}_s^{t,\zeta}|^2\dif r+C\left(\sup\limits_{s\in[t,T]}\mE\sup _{s \leq r \leq s+\d}|Y_r^{\e,t,\zeta}-Y_{s+\d}^{\e,t,\zeta}|^{2}\right)\no\\
&&+2L^2_4TC(1+\mE|\zeta|^{4q_3})^{1/2}(\d+\d^2).
\label{i2}
\ee

For $I_3$, by the same deduction to that for $I_2$, we have that
\be
I_3&\leq&\int_t^T\sup\limits_{s\in[t,r]}\mE|Y_s^{\e,t,\zeta}-\bar{Y}_s^{t,\zeta}|^2\dif r+C\left(\sup\limits_{s\in[t,T]}\mE\sup _{s \leq r \leq s+\d}|Y_r^{\e,t,\zeta}-Y_{s+\d}^{\e,t,\zeta}|^{2}\right)\no\\
&&+2L^2_4TC(1+\mE|\zeta|^{4q_3})^{1/2}(\d+\d^2).
\label{i3}
\ee

For $I_4$, we define $X^{\e,t,\zeta}_r=\zeta, Y^{\e,t,\zeta}_r=Y^{\e,t,\zeta}_t, \bar{Y}_{r}^{t,\zeta}=\bar{Y}_{t}^{t,\zeta}$ for $r\in[0,t]$ and obtain that
\ce
I_4&\leq&2\mE\left|\int_0^T\<Y_{r(\d)}^{\e,t,\zeta}-\bar{Y}_{r(\d)}^{t,\zeta},f_1(\frac{r}{\e},X^{\e,t,\zeta}_{r(\d)},Y_{r(\d)}^{\e,t,\zeta})-\bar{f}_1(X^{\e,t,\zeta}_{r(\d)},Y_{r(\d)}^{\e,t,\zeta})\>\dif r\right|\\
&&+2\sup\limits_{s\in[0,T]}\mE\left|\int_0^s\<Y_{r(\d)}^{\e,t,\zeta}-\bar{Y}_{r(\d)}^{t,\zeta},f_1(\frac{r}{\e},X^{\e,t,\zeta}_{r(\d)},Y_{r(\d)}^{\e,t,\zeta})-\bar{f}_1(X^{\e,t,\zeta}_{r(\d)},Y_{r(\d)}^{\e,t,\zeta})\>\dif r\right|\\
&\leq&4\sup\limits_{s\in[0,T]}\mE\left|\int_0^s\<Y_{r(\d)}^{\e,t,\zeta}-\bar{Y}_{r(\d)}^{t,\zeta},f_1(\frac{r}{\e},X^{\e,t,\zeta}_{r(\d)},Y_{r(\d)}^{\e,t,\zeta})-\bar{f}_1(X^{\e,t,\zeta}_{r(\d)},Y_{r(\d)}^{\e,t,\zeta})\>\dif r\right|\\
&\leq&4\sup\limits_{s\in[0,T]}\mE\left|\int_0^{[\frac{s}{\d}]\d}\<Y_{r(\d)}^{\e,t,\zeta}-\bar{Y}_{r(\d)}^{t,\zeta},f_1(\frac{r}{\e},X^{\e,t,\zeta}_{r(\d)},Y_{r(\d)}^{\e,t,\zeta})-\bar{f}_1(X^{\e,t,\zeta}_{r(\d)},Y_{r(\d)}^{\e,t,\zeta})\>\dif r\right|\\
&&+4\sup\limits_{s\in[0,T]}\mE\left|\int_{[\frac{s}{\d}]\d}^s\<Y_{r(\d)}^{\e,t,\zeta}-\bar{Y}_{r(\d)}^{t,\zeta},f_1(\frac{r}{\e},X^{\e,t,\zeta}_{r(\d)},Y_{r(\d)}^{\e,t,\zeta})-\bar{f}_1(X^{\e,t,\zeta}_{r(\d)},Y_{r(\d)}^{\e,t,\zeta})\>\dif r\right|\\
&=:&I_{41}+I_{42}.
\de
Then $(\mathbf{H}^{1}_{f})$ yields that
\be
I_{41}&\leq&4\sup\limits_{s\in[0,T]}\mE\sum\limits_{k=0}^{[\frac{s}{\d}]-1}\left|\int_{k\d}^{(k+1)\d}\<Y_{k\d}^{\e,t,\zeta}-\bar{Y}_{k\d}^{t,\zeta},f_1(\frac{r}{\e},X^{\e,t,\zeta}_{k\d},Y_{k\d}^{\e,t,\zeta})-\bar{f}_1(X^{\e,t,\zeta}_{k\d},Y^{\e,t,\zeta}_{k\d})\>\dif r\right|\no\\
&\leq&4\sum\limits_{k=0}^{[\frac{T}{\d}]-1}\mE\left|\int_{k\d}^{(k+1)\d}\<Y_{k\d}^{\e,t,\zeta}-\bar{Y}_{k\d}^{t,\zeta},f_1(\frac{r}{\e},X^{\e,t,\zeta}_{k\d},Y_{k\d}^{\e,t,\zeta})-\bar{f}_1(X^{\e,t,\zeta}_{k\d},Y_{k\d}^{\e,t,\zeta})\>\dif r\right|\no\\
&\leq&4[\frac{T}{\d}]\sup\limits_{0\leq k\leq [\frac{T}{\d}]-1}\left(\mE|Y_{k\d}^{\e,t,\zeta}-\bar{Y}_{k\d}^{t,\zeta}|^2\right)^{1/2}\no\\
&&\times\left(\mE\left|\int_{k\d}^{(k+1)\d}\left(f_1(\frac{r}{\e},X^{\e,t,\zeta}_{k\d},Y_{k\d}^{\e,t,\zeta})-\bar{f}_1(X^{\e,t,\zeta}_{k\d},Y_{k\d}^{\e,t,\zeta})\right)\dif r\right|^2\right)^{1/2}\no\\
&\leq&4C[\frac{T}{\d}]\d\sup\limits_{0\leq k\leq [\frac{T}{\d}]-1}\left(\mE\left|\frac{\e}{\d}\int_{k\d/\e}^{(k+1)\d/\e}\left(f_1(u,X^{\e,t,\zeta}_{k\d},Y_{k\d}^{\e,t,\zeta})-\bar{f}_1(X^{\e,t,\zeta}_{k\d},Y_{k\d}^{\e,t,\zeta})\right)\dif u\right|^2\right)^{1/2}\no\\
&\leq&4CT\kappa_3^{1/2}(\frac{\d}{\e}).
\label{i41}
\ee
And by the H\"older inequality, it holds that
\be
I_{42}&\leq& 4\sup\limits_{s\in[0,T]}\left(\int_{[\frac{s}{\d}]\d}^s\mE|f_1(\frac{r}{\e},X^{\e,t,\zeta}_{r(\d)},Y_{r(\d)}^{\e,t,\zeta})-\bar{f}_1(X^{\e,t,\zeta}_{r(\d)},Y_{r(\d)}^{\e,t,\zeta})|^2\dif r\right)^{1/2}\no\\
&&\times\left(\int_{[\frac{s}{\d}]\d}^s\mE|Y_{r(\d)}^{\e,t,\zeta}-\bar{Y}_{r(\d)}^{t,\zeta}|^2\dif r\right)^{1/2}\no\\
&\leq& 4C\sup\limits_{s\in[0,T]}\left(\int_{[\frac{s}{\d}]\d}^s(1+\mE|X^{\e,t,\zeta}_{r(\d)}|^{2q_3+2}+\mE|Y_{r(\d)}^{\e,t,\zeta})|^2)\dif r\right)^{1/2}\no\\
&&\times\left(\int_{[\frac{s}{\d}]\d}^s(\mE|Y_{r(\d)}^{\e,t,\zeta}|^2+\mE|\bar{Y}_{r(\d)}^{t,\zeta}|^2)\dif r\right)^{1/2}\no\\
&\leq&C\d.
\label{i42}
\ee
Combining (\ref{i1})-(\ref{i42}) with (\ref{i1234}), we conclude (\ref{yyfbarf}). The proof is complete.

\section{Applications}\label{app}

In this section, we apply Theorem \ref{averprin} to nonlinear parabolic PDEs.

\subsection{Application to the viscosity solutions of multivalued PDEs}

In this subsection, we require that $d=1, \zeta=x\in\mR^m$ and study the homogenization for multivalued PDEs.

Let $u^\e$ be the viscosity solution of the PDE
\be\left\{\begin{array}{l}
\frac{\p u^\e(t,x)}{\p t}+\sL^\e u^\e(t,x)+f_1(\frac{t}{\e},x,u^\e(t,x))+f_2(\nabla u^\e(t,x)\s(\frac{t}{\e},x))\in \p \varphi(u^\e(t,x)),~t\in[0,T], \\ 
u^\e(T,x)=g(x), \quad u^\e(t,x)\in\overline{{\rm Dom}(\varphi)}, \quad x\in\mR^m,
\end{array}
\right.
\label{1dpde1}
\ee
where 
\ce
\sL^\e:=\frac{1}{2}\sum\limits_{i,j=1}^m (\s\s^*)_{ij}(\frac{t}{\e},x)\frac{\p^2}{\p x_i \p x_j}+\sum\limits_{i=1}^m b_i(\frac{t}{\e},x)\frac{\p}{\p x_i},
\de
and $\bar{u}$ be the viscosity solution of the PDE
\be\left\{\begin{array}{l}
\frac{\p \bar{u}(t,x)}{\p t}+\bar{\sL} \bar{u}(t,x)+\bar{f}_1(x,\bar{u}(t,x))+f_2(\nabla \bar{u}(t,x)\bar{\s}(x))\in \p \varphi(\bar{u}(t,x)), \quad t\in[0,T], \\ 
\bar{u}(T,x)=g(x), \quad \bar{u}(t,x)\in\overline{{\rm Dom}(\varphi)}, \quad   x\in\mR^m,
\end{array}
\right.
\label{1dpde2}
\ee
where 
\ce
\bar{\sL}:=\frac{1}{2}\sum\limits_{i,j=1}^m (\bar{\s}\bar{\s}^*)_{ij}(x)\frac{\p^2}{\p x_i \p x_j}+\sum\limits_{i=1}^m \bar{b}_i(x)\frac{\p}{\p x_i}.
\de

Assume:
\begin{enumerate}[$(\mathbf{H}^3_{f})$]
\item For each $R>0$ there exists a continuous function $\a_R: \mR_+\rightarrow \mR_+, \a_R(0)=0$ such that 
\ce
&&|f_1(s,x_1,y)+f_2(z)-f_1(s,x_2,y)-f_2(z)|\leq \a_R((|x_1-x_2|)(1+|z|)), \\
&&\quad s\in[0,T], \quad  |x_1|, |x_2|\leq R, \quad z\in\mR^l.
\de
\end{enumerate}

\bt\label{hoth1}
Assume that $(\mathbf{H}^{1}_{b, \s})$, $(\mathbf{H}_{\varphi})$, $(\mathbf{H}_{g})$, $(\mathbf{H}^{1}_{f})$, $(\mathbf{H}^{2}_{b, \s})$, $(\mathbf{H}^{2}_{f})$ and $(\mathbf{H}^3_{f})$ hold. Then for any $(t,x)\in[0,T]\times\mR^m$
$$
\lim\limits_{\e\rightarrow 0}u^\e(t,x)=\bar{u}(t,x).
$$
\et
\begin{proof}
First of all, by \cite[Theorem 4.1 and 4.2]{pr}, we know that $u^\e(t,x)=Y_t^{\e,t,x}, \bar{u}(t,x)=\bar{Y}_t^{t,x}$ are unique viscosity solutions of Eq.(\ref{1dpde1}) and Eq.(\ref{1dpde2}), respectively. Therefore, the result follows from Theorem \ref{averprin}.
\end{proof}

\br
If $f_1(s,x,y)$ is independent of $s$ and $f_2(z)=0$, Theorem \ref{hoth1} is the same to \cite[Theorem 4.1]{eo}. Therefore, our result is more general.
\er

\subsection{Application to the viscosity solutions for systems of multivalued PDEs}

In this subsection, we require that $\zeta=x\in\mR^m, f_2(z)=0$ and study the homogenization for systems of multivalued PDEs.

Let $v^\e$ be the viscosity solution of the PDE
\be\left\{\begin{array}{l}
\frac{\p v^\e(t,x)}{\p t}+\sL^\e v^\e(t,x)+f_1(\frac{t}{\e},x,v^\e(t,x))\in \p \varphi(v^\e(t,x)), \quad t\in[0,T], \\ 
v^\e(T,x)=g(x), \quad v^\e(t,x)\in\overline{{\rm Dom}(\varphi)}, \quad x\in\mR^m,
\end{array}
\right.
\label{ddpde1}
\ee
where 
\ce
(\sL^\e v_k^\e)(t,x):=\frac{1}{2}\sum\limits_{i,j=1}^m (\s\s^*)_{ij}(\frac{t}{\e},x)\frac{\p^2 v_k^\e(t,x)}{\p x_i \p x_j}+\sum\limits_{i=1}^m b_i(\frac{t}{\e},x)\frac{\p v_k^\e(t,x)}{\p x_i}, \quad k=1,\cdots,d,
\de
and $\bar{v}$ be the viscosity solution of the PDE
\be\left\{\begin{array}{l}
\frac{\p \bar{v}(t,x)}{\p t}+\bar{\sL} \bar{v}(t,x)+\bar{f}_1(x,\bar{v}(t,x))\in \p \varphi(\bar{v}(t,x)), \quad t\in[0,T], \\ 
\bar{v}(T,x)=g(x), \quad \bar{v}(t,x)\in\overline{{\rm Dom}(\varphi)},  \quad x\in\mR^m,
\end{array}
\right.
\label{ddpde2}
\ee
where 
\ce
\bar{\sL}\bar{v}_k(t,x):=\frac{1}{2}\sum\limits_{i,j=1}^m (\bar{\s}\bar{\s}^*)_{ij}(x)\frac{\p^2}{\p x_i \p x_j}+\sum\limits_{i=1}^m \bar{b}_i(x)\frac{\p}{\p x_i}, \quad k=1,\cdots,d.
\de

Assume:
\begin{enumerate}[$(\mathbf{H}^{\prime}_{\varphi})$]
\item For all $\rho \in Dom(\varphi)$, there exists a neighborhood $V$ of $\rho$ such that $(\partial \varphi)_0$ is bounded on $\cD(\partial \varphi) \cap V$, where for $\varrho\in \cD(\p\varphi)$, $(\partial \varphi)_0(\varrho)\in\mR^d$ is the unique one such that $|(\partial \varphi)_0(\varrho)|=\inf|\p\varphi(\varrho)|$. Moreover, if $\rho \in Dom(\varphi)$ and $y \in \mathbb{R}^d$ such that $\rho+y \in Dom(\varphi)$, there exists a neighbourhood $V$ of $\rho$ such that
$$
\forall \vartheta \in V \cap \cD(\partial \varphi), \exists t>0: \quad \vartheta+t y \in \cD(\partial \varphi).
$$
\end{enumerate}

\bt\label{hoth2}
Assume that $(\mathbf{H}^{1}_{b, \s})$, $(\mathbf{H}_{g})$, $(\mathbf{H}^{1}_{f})$, $(\mathbf{H}^{2}_{b, \s})$, $(\mathbf{H}^{2}_{f})$ and $(\mathbf{H}^{\prime}_{\varphi})$ hold. Then for any $(t,x)\in[0,T]\times\mR^m$
$$
\lim\limits_{\e\rightarrow 0}u^\e(t,x)=\bar{u}(t,x).
$$
\et
\begin{proof}
First of all, \cite[Theorem 6 and 7]{mprz} imply that $v^\e(t,x)=Y_t^{\e,t,x}, \bar{v}(t,x)=\bar{Y}_t^{t,x}$ are unique viscosity solutions of Eq.(\ref{ddpde1}) and Eq.(\ref{ddpde2}), respectively. Therefore, the result follows from Theorem \ref{averprin}.
\end{proof}

\subsection{Application to the viscosity solutions for quasilinear parabolic PDEs}

In this subsection, we require that $\varphi=0, d=1, \zeta=x\in\mR^m$ and study the homogenization for quasilinear parabolic PDEs.

Let $w^\e$ be the viscosity solution of the PDE
\be\left\{\begin{array}{l}
\frac{\p w^\e(t,x)}{\p t}+\sL^\e w^\e(t,x)+f_1(\frac{t}{\e},x,w^\e(t,x))+f_2(\nabla w^\e(t,x)\s(\frac{t}{\e},x))=0, \quad t\in[0,T], \\ 
w^\e(T,x)=g(x), \quad x\in\mR^m,
\end{array}
\right.
\label{1dpdephi1}
\ee
and $\bar{w}$ be the viscosity solution of the PDE
\be\left\{\begin{array}{l}
\frac{\p \bar{w}(t,x)}{\p t}+\bar{\sL} \bar{w}(t,x)+\bar{f}_1(x,\bar{w}(t,x))+f_2(\nabla \bar{w}(t,x)\bar{\s}(x))=0, \quad t\in[0,T], \\ 
\bar{w}(T,x)=g(x), \quad x\in\mR^m.
\end{array}
\right.
\label{1dpdephi2}
\ee

By \cite[Theorem 4.3]{pp2}, we conclude the relationship between $w^\e$ and $\bar{w}$ as follows.

\bt\label{hoth2}
Assume that $(\mathbf{H}^{1}_{b, \s})$, $(\mathbf{H}_{g})$, $(\mathbf{H}^{1}_{f})$, $(\mathbf{H}^{2}_{b, \s})$  and $(\mathbf{H}^{2}_{f})$ hold. Then for any $(t,x)\in[0,T]\times\mR^m$
$$
\lim\limits_{\e\rightarrow 0}u^\e(t,x)=\bar{u}(t,x).
$$
\et

\section{An example}\label{exam}

In this section, we give an example to illustrate our result.

\bx\label{ex1}
Assume that $m=l=d=1, \zeta=x\in\mR$ and consider the following forward-backward multivalued stochastic system: $0\leq t\leq s\leq T$
\be\left\{\begin{array}{l}
\dif X_{s}^{\e,t,x}=\frac{\frac{s}{\e}}{1+\frac{s}{\e}}\cos(X_{s}^{\e,t,x})\dif s+(1-e^{-\frac{s}{2\e}})\sin(X_{s}^{\e,t,x})\dif W_{s},\\
X_{t}^{\e,t,x}=x,\\
\dif Y_{s}^{\e,t,x}\in \p\varphi(Y_{s}^{\e,t,x})\dif s+Z_{s}^{\e,t,x}\dif W_{s},\\
Y_{T}^{\e,t,x}=g(X_{T}^{\e,t,x})\in\overline{\cD(\p\varphi)},
\end{array}
\right.
\label{multfbsdeex1}
\ee
where $\varphi, g$ satisfy $(\mathbf{H}_{\varphi})$, $(\mathbf{H}_{g})$. It is easy to see that $b(s,x)=\frac{s}{1+s}\cos(x), \s(s,x)=(1-e^{-s/2})\sin(x), f_1(s,x,y)=f_2(z)=0$ satisfy $(\mathbf{H}^{1}_{b, \s})$ and $(\mathbf{H}^{1}_{f})$. Thus, the system (\ref{multfbsdeex1}) has a unique solution $(X^{\e,t,x}, Y^{\e,t,x}, K^{\e,t,x}, Z^{\e,t,x})$.

Taking $\bar{b}(x)=\cos(x), \bar{\s}(x)=\sin(x)$, we justify that
\ce
&&\left|\frac{1}{\hat{T}}\int_0^{\hat{T}}b(s,x)\dif s-\bar{b}(x)\right|^2=\left|\frac{1}{\hat{T}}\int_0^{\hat{T}}\frac{s}{1+s}\cos(x)\dif s-\cos(x)\right|^2\\
&\leq&\frac{1}{\hat{T}}\int_0^{\hat{T}}\frac{1}{(1+s)^2}\dif s|x|^2\leq\frac{1}{\hat{T}}(1+|x|^2),
\de
and
\ce
&&\frac{1}{\hat{T}}\int_0^{\hat{T}}\left|\s(s,x)-\bar{\s}(x)\right|^2\dif s=\frac{1}{\hat{T}}\int_0^{\hat{T}}\left|(1-e^{-s/2})\sin(x)-\sin(x)\right|^2\dif s\\
&\leq&\frac{1}{\hat{T}}\int_0^{\hat{T}}e^{-s}\dif s|x|^2\leq\frac{1}{\hat{T}}(1+|x|^2).
\de
That is, $(\mathbf{H}^{2}_{b, \s})$ holds. Therefore, by Theorem \ref{averprin}, we obtain that 
\ce
\lim\limits_{\e\rightarrow0}\mE\sup\limits_{s\in[t,T]}|Y^{\e,t,x}_s-\bar{Y}^{t,x}_{s}|^2=0,
\de
where $(\bar{X}^{t,x}, \bar{Y}^{t,x}, \bar{K}^{t,x}, \bar{Z}^{t,x})$ is the unique solution of the following system:
\ce\left\{\begin{array}{l}
\dif \bar{X}_{s}^{t,x}=\cos(\bar{X}_{s}^{t,x})\dif s+\sin(\bar{X}_{s}^{t,x})\dif W_{s},\\
\bar{X}_{t}^{t,x}=x,\\
\dif \bar{Y}_{s}^{t,x}\in \p\varphi(\bar{Y}_{s}^{t,x})\dif s+\bar{Z}_{s}^{t,x}\dif W_{s},\\
\bar{Y}_{T}^{t,x}=g(\bar{X}_{T}^{t,x})\in\overline{\cD(\p\varphi)}.
\end{array}
\right.
\de

Next, we investigate the following PDEs:
\be\left\{\begin{array}{l}
\frac{\p u^\e(t,x)}{\p t}+\frac{1}{2}(1-e^{-\frac{t}{2\e}})^2\sin^2(x)\frac{\p^2 u^\e(t,x)}{\p x \p x}+\frac{\frac{t}{\e}}{1+\frac{t}{\e}}\cos(x)\frac{\p u^\e(t,x)}{\p x}\in \p \varphi(u^\e(t,x)), t\in[0,T], \\ 
u^\e(T,x)=g(x), \quad x\in\mR,
\end{array}
\right.
\label{pde1}
\ee
and 
\be\left\{\begin{array}{l}
\frac{\p \bar{u}(t,x)}{\p t}+\frac{1}{2}\sin^2(x)\frac{\p^2 \bar{u}(t,x)}{\p x \p x}+\cos(x)\frac{\p \bar{u}(t,x)}{\p x}\in \p \varphi(\bar{u}(t,x)), t\in[0,T], \\ 
\bar{u}(T,x)=g(x), \quad x\in\mR.
\end{array}
\right.
\label{pde2}
\ee
By Theorem \ref{hoth1}, it holds that the viscosity solution $u^\e(t,x)$ of Eq.(\ref{pde1}) converges to the viscosity solution $\bar{u}(t,x)$ of Eq.(\ref{pde2}).
\ex

\end{document}